\documentclass[11pt,a4paper, reqno]{amsart}
\usepackage[utf8]{inputenc}       
\usepackage[
a4paper,
left=4cm,
right=4cm,
top=2.9cm,
bottom=2.9cm
]{geometry}
\usepackage{natbib}
\usepackage{csquotes}

\usepackage{amsmath, amsthm, amssymb}
\usepackage{setspace}
\usepackage{tikz-cd}
\usetikzlibrary{
	knots,
	hobby,
	decorations.pathreplacing,
	shapes.geometric,
	calc,
	decorations.markings
}

\usepackage{etoolbox}

\makeatletter
\let\old@tocline\@tocline
\let\section@tocline\@tocline

\newcommand{\subsection@dotsep}{4.5}
\newcommand{\subsubsection@dotsep}{4.5}
\patchcmd{\@tocline}
{\hfil}
{\nobreak
	\leaders\hbox{$\m@th
		\mkern \subsection@dotsep mu\hbox{.}\mkern \subsection@dotsep mu$}\hfill
	\nobreak}{}{}
\let\subsection@tocline\@tocline
\let\@tocline\old@tocline

\patchcmd{\@tocline}
{\hfil}
{\nobreak
	\leaders\hbox{$\m@th
		\mkern \subsubsection@dotsep mu\hbox{.}\mkern \subsubsection@dotsep mu$}\hfill
	\nobreak}{}{}
\let\subsubsection@tocline\@tocline
\let\@tocline\old@tocline

\let\old@l@subsection\l@subsection
\let\old@l@subsubsection\l@subsubsection

\def\@tocwriteb#1#2#3{%
	\begingroup
	\@xp\def\csname #2@tocline\endcsname##1##2##3##4##5##6{%
		\ifnum##1>\c@tocdepth
		\else \sbox\z@{##5\let\indentlabel\@tochangmeasure##6}\fi}%
	\csname l@#2\endcsname{#1{\csname#2name\endcsname}{\@secnumber}{}}%
	\endgroup
	\addcontentsline{toc}{#2}%
	{\protect#1{\csname#2name\endcsname}{\@secnumber}{#3}}}%

\newlength{\@tocsectionindent}
\newlength{\@tocsubsectionindent}
\newlength{\@tocsubsubsectionindent}
\newlength{\@tocsectionnumwidth}
\newlength{\@tocsubsectionnumwidth}
\newlength{\@tocsubsubsectionnumwidth}
\newcommand{\settocsectionnumwidth}[1]{\setlength{\@tocsectionnumwidth}{#1}}
\newcommand{\settocsubsectionnumwidth}[1]{\setlength{\@tocsubsectionnumwidth}{#1}}
\newcommand{\settocsubsubsectionnumwidth}[1]{\setlength{\@tocsubsubsectionnumwidth}{#1}}
\newcommand{\settocsectionindent}[1]{\setlength{\@tocsectionindent}{#1}}
\newcommand{\settocsubsectionindent}[1]{\setlength{\@tocsubsectionindent}{#1}}
\newcommand{\settocsubsubsectionindent}[1]{\setlength{\@tocsubsubsectionindent}{#1}}

\renewcommand{\l@section}{\section@tocline{1}{\@tocsectionvskip}{\@tocsectionindent}{\@tocsectionnumwidth}{\@tocsectionformat}}%
\renewcommand{\l@subsection}{\subsection@tocline{1}{\@tocsubsectionvskip}{\@tocsubsectionindent}{\@tocsubsectionnumwidth}{\@tocsubsectionformat}}%
\renewcommand{\l@subsubsection}{\subsubsection@tocline{1}{\@tocsubsubsectionvskip}{\@tocsubsubsectionindent}{\@tocsubsubsectionnumwidth}{\@tocsubsubsectionformat}}%
\newcommand{\@tocsectionformat}{}
\newcommand{\@tocsubsectionformat}{}
\newcommand{\@tocsubsubsectionformat}{}
\expandafter\def\csname toc@1format\endcsname{\@tocsectionformat}
\expandafter\def\csname toc@2format\endcsname{\@tocsubsectionformat}
\expandafter\def\csname toc@3format\endcsname{\@tocsubsubsectionformat}
\newcommand{\settocsectionformat}[1]{\renewcommand{\@tocsectionformat}{#1}}
\newcommand{\settocsubsectionformat}[1]{\renewcommand{\@tocsubsectionformat}{#1}}
\newcommand{\settocsubsubsectionformat}[1]{\renewcommand{\@tocsubsubsectionformat}{#1}}
\newlength{\@tocsectionvskip}
\newcommand{\settocsectionvskip}[1]{\setlength{\@tocsectionvskip}{#1}}
\newlength{\@tocsubsectionvskip}
\newcommand{\settocsubsectionvskip}[1]{\setlength{\@tocsubsectionvskip}{#1}}
\newlength{\@tocsubsubsectionvskip}
\newcommand{\settocsubsubsectionvskip}[1]{\setlength{\@tocsubsubsectionvskip}{#1}}

\patchcmd{\tocsection}{\indentlabel}{\makebox[\@tocsectionnumwidth][l]}{}{}
\patchcmd{\tocsubsection}{\indentlabel}{\makebox[\@tocsubsectionnumwidth][l]}{}{}
\patchcmd{\tocsubsubsection}{\indentlabel}{\makebox[\@tocsubsubsectionnumwidth][l]}{}{}

\newcommand{\@sectypepnumformat}{}
\renewcommand{\contentsline}[1]{%
	\expandafter\let\expandafter\@sectypepnumformat\csname @toc#1pnumformat\endcsname%
	\csname l@#1\endcsname}
\newcommand{\@tocsectionpnumformat}{}
\newcommand{\@tocsubsectionpnumformat}{}
\newcommand{\@tocsubsubsectionpnumformat}{}
\newcommand{\setsectionpnumformat}[1]{\renewcommand{\@tocsectionpnumformat}{#1}}
\newcommand{\setsubsectionpnumformat}[1]{\renewcommand{\@tocsubsectionpnumformat}{#1}}
\newcommand{\setsubsubsectionpnumformat}[1]{\renewcommand{\@tocsubsubsectionpnumformat}{#1}}
\renewcommand{\@tocpagenum}[1]{%
	\hfill {\mdseries\@sectypepnumformat #1}}

\let\oldappendix\appendix
\renewcommand{\appendix}{%
	\leavevmode\oldappendix%
	\addtocontents{toc}{%
		\protect\settowidth{\protect\@tocsectionnumwidth}{\protect\@tocsectionformat\sectionname\space}%
		\protect\addtolength{\protect\@tocsectionnumwidth}{2em}}%
}
\makeatother

\makeatletter
\settocsectionnumwidth{2em}
\settocsubsectionnumwidth{2.5em}
\settocsubsubsectionnumwidth{3em}
\settocsectionindent{1pc}%
\settocsubsectionindent{\dimexpr\@tocsectionindent+\@tocsectionnumwidth}%
\settocsubsubsectionindent{\dimexpr\@tocsubsectionindent+\@tocsubsectionnumwidth}%
\makeatother

\settocsectionvskip{10pt}
\settocsubsectionvskip{0pt}
\settocsubsubsectionvskip{0pt}

\settocsectionformat{\bfseries}%
\settocsubsectionformat{\mdseries}
\settocsubsubsectionformat{\mdseries}
\setsectionpnumformat{}
\setsubsectionpnumformat{\mdseries}
\setsubsubsectionpnumformat{\mdseries}

\theoremstyle{plain}
\theoremstyle{definition}
\newtheorem{Theorem}{Theorem}[section]
\definecolor{shadecolor}{RGB}{217, 255, 217}

\newtheorem{Lemma}[Theorem]{Lemma}

\newtheorem{Proposition}[Theorem]{Proposition}

\newtheorem{Corollary}[Theorem]{Corollary}

\numberwithin{equation}{section}

\DeclareMathOperator{\im}{Im}
\DeclareMathOperator{\supp}{supp}
\DeclareMathOperator{\cl}{clos}

\DeclareMathOperator{\sinc}{sinc}
\newcommand{\pw}{\mathrm{PW}_a^p}
\newcommand{\sr}{\mathcal{S}(\R)}
\newcommand{\sar}{\mathcal{S}_a(\R)}

\usepackage[colorlinks=blue, urlcolor=blue, citecolor=black]{hyperref}

\renewcommand{\P}{\mathbb{P}}

\newcommand{\D}{\mathbb{D}}

\newcommand{\C}{\mathbb{C}}
\renewcommand{\S}{\mathcal{S}}

\newcommand{\T}{\mathbb{T}}
\newcommand{\F}{\mathcal{F}}
\newcommand{\R}{\mathbb{R}}
\newcommand{\Z}{\mathbb{Z}}
\newcommand{\N}{\mathbb{N}}

\DeclareFontFamily{U}{mathx}{}
\DeclareFontShape{U}{mathx}{m}{n}{<-> mathx10}{}
\DeclareSymbolFont{mathx}{U}{mathx}{m}{n}
\DeclareMathAccent{\widehat}{0}{mathx}{"70}
\DeclareMathAccent{\widecheck}{0}{mathx}{"71}

\title[Bounded symbols of Toeplitz operators and a weak factorization]{Bounded symbols of Toeplitz operators on Paley-Wiener spaces and a weak factorization theorem}
\author[Petr Kulikov]{Petr Kulikov}

\begin{document}
	\thispagestyle{empty}
	\markboth{\scshape PETR KULIKOV}{\scshape BOUNDED SYMBOLS OF TOEPLITZ OPERATORS, WEAK FACTORIZATION}
	
	\begin{center}
		{\sc\bf BOUNDED SYMBOLS OF TOEPLITZ OPERATORS ON PALEY-WIENER SPACES AND A WEAK FACTORIZATION THEOREM}
	\end{center}
	
	\vspace{0.4cm}
		
	\begin{center}
		{\sc\small PETR KULIKOV}
	\end{center}
	
	\vspace{0.4cm}
	
	\begin{quote}\fontsize{9pt}{11pt}\selectfont
		{\sc Abstract.} A classical result by R. Rochberg says that every bounded Toeplitz operator $T$ on the Hilbert Paley-Wiener space $\mathrm{PW}_a^2$ admits a bounded symbol $\varphi$. We generalize this result to Toeplitz operators on the Banach Paley-Wiener spaces $\mathrm{PW}_a^p$, $1<p<+\infty$. 
				
		The Toeplitz commutator theorem describes the integral identity that must hold for a bounded operator $T$ on $\pw$ to be a Toeplitz operator on $\pw$. We prove this theorem in the continuous case, thus extending the result previously obtained by D. Sarason in the discrete case.
		
		Upon combining the results, we establish the weak factorization theorem, namely, for $p,q>1$, $\frac{1}{p}+\frac{1}{q}=1$, any function $h$ belonging to $\mathrm{PW}^1_{2a}$ can be represented as $$h=\sum_{k\geqslant 0}f_k\bar{g}_k,\qquad f_k\in\pw,\,g_k\in\mathrm{PW}_a^q.$$
	\end{quote}
	
	\vspace{0.4cm}
	
	
	\medskip

	\section{Introduction}\label{ps}
	Let $\sr$ denote the classical Schwartz space. The Fourier transform on $\sr$ is denoted by $\F[f](\xi)=\widehat{f}(\xi)$. Fix $a > 0$ and define the set 
	$$\sar=\{f\in \sr\mid \supp{\widehat{f}}\subset[-a,a]\}.$$ 
	For $1\leqslant p <+\infty$, the Paley-Wiener space $\pw$ is a closed subspace of $L^p(\R)$ defined by $\pw = \cl_{L^p(\R)}\sar$. In particular,
	$$\mathrm{PW}_a^2 = \{f \in L^2(\R)\mid \widehat{f}=0\text{ a.e.\;on }\R\setminus [-a,a]\}.$$
	Let $\mathfrak{m}$ be a bounded measurable function on $\R$. The Fourier multiplier associated to symbol $\mathfrak{m}$ is the map defined by $f \mapsto \F^{-1}\mathfrak{m} \F [f]$ for any $f\in\sr$. Fix the Fourier multiplier associated to symbol $\chi_{[-a,a]}$ and denote it by $\P_a$,
	$$\P_a \;:\; f\longmapsto \F^{-1}\chi_{[-a,a]} \F [f],\qquad f\in \sr. $$
	Since $\chi_{[-a,a]}^2 = \chi_{[-a,a]}$, $\P_a$ is, in fact, a linear bounded projector to $\pw$.	
	
	\medskip
	
	Let $\mathcal{P}(\R)$ denote the set of all complex-valued functions defined on $\R$ that grow no faster than polynomials:
	$$\mathcal{P}(\R) = \{f \colon\R\to\C \mid \exists n \in \N:  \sup\limits_{x\in \R}|f(x)|\cdot(1 + |x|)^{-n} <  +\infty\}.$$
	A Toeplitz operator $T_\varphi \,\colon \, \pw \to \pw$ with symbol $\varphi \in \mathcal{P}(\R)$ is a mapping densely defined by 
	$$T_\varphi\; \colon\; f  \mapsto \P_a[\varphi \cdot f],\qquad f\in\sar .$$
	Since $\mathcal{P}(\R) \cdot \sar \subset L^p(\R)$, we have $\varphi \cdot f \in L^p(\R)$ for every $f\in\sar$. Hence, $T_\varphi$ is well defined. In the case
	$$\sup \{\left \|T_\varphi [f]\right \|_{L^p(\R)}\mid f\in\sar, \left \|f\right\|_{L^p(\R)}=1\} <+\infty,$$
	the operator $T_\varphi$ admits a unique bounded extension to $\pw$. This extension will be denoted by the same notation $T_\varphi$.
	
	\medskip
	
	The symbol of a Toeplitz operator on $\pw$ is not unique. We say that a Toeplitz operator $T_\varphi$ on $\pw$ admits a bounded symbol $\psi$ if $T_\varphi=T_\psi$ for a function $\psi \in L^\infty(\R)$. Clearly, any bounded symbol $\varphi \in L^\infty(\R)$ determines the bounded Toeplitz operator $T_\varphi$ on $\pw$, and
	$$\left \|T_\varphi \right \|_{\pw\to\pw} \leqslant \left \| \varphi \right \|_{L^\infty(\R)}.$$
	The class of all bounded Toeplitz operators on $\pw$ will be denoted by $\mathcal{T}^p(a)$. It is easy to see that some unbounded symbols $\varphi$ can produce bounded Toeplitz operators on $\pw$. For instance, this is the case for the symbol 
	$$\varphi(x) = x\cdot e^{2\pi i x \cdot 2a}, \quad x\in\R.$$ 
	Indeed, for every $f \in \sar$ we have $\supp \F[\varphi\cdot f]\subset [a, 3a]$. Thus, $\P_a[\varphi\cdot f] = 0$ and $T_\varphi = 0$ as an operator on $\pw$. This motivates the question of whether every bounded Toeplitz operator on $\pw$ admits a bounded symbol. In the case $p=2$, the affirmative answer to this question was given by R. Rochberg \citep{RR} in 1987.
	
	\medskip
	
	Our aim in the present paper is to prove the following theorems.
	
	\begin{Theorem}\label{result} \textit{
		Let $1<p<+\infty$. Let $T_\varphi$ be a Toeplitz operator on $\pw$ with symbol $\varphi \in \sr$. Then $T_\varphi$ admits a bounded symbol $\psi$ such that 
		$$\left \| \psi \right \|_{L^\infty(\R)} \leqslant c \biggl(p+\frac{1}{p-1}\biggr) \cdot  \left \| T_\varphi \right \|_{\pw\to\pw},$$
		for a universal constant $c > 0$.} 
	\end{Theorem}
	
	\begin{Theorem}\label{result2} \textit{
		Let $1 < p < \infty$ and $\tfrac{1}{p}+\tfrac{1}{q}=1$. For any function $h\in \mathrm{PW}_{2a}^1$ there exist $f_k \in \mathrm{PW}_a^p$, $g_k \in \mathrm{PW}_a^q$ with $$\sum\limits_{k=0}^\infty \left \|f_k\right \|_{L^p(\R)}\cdot\left \|g_k\right \|_{L^q(\R)}<+\infty \;\text{ such that }\;h = \sum\limits_{k=0}^\infty f_k \,\overline{g}_k.$$}
	\end{Theorem}
	A result similar to Theorem \ref{result2} is usually called a weak factorization \text{theorem.} 
	
	\medskip
	
	\subsection{Notations}
	We normalize the Fourier transform on $\sr$ by
	$$\F[f](\xi)=\widehat{f}(\xi) = \int\limits_\R  e^{-2\pi i \xi x} f(x) \,dx, \qquad \xi\in\R.$$
	For the inverse Fourier transform associated with the one defined above we use the following notion: $\mathcal{F}^{-1}[f] = \widecheck{f}$.
	
	For $f_1,f_2\in L^2(\R)$, the dual pairing is given by
	$$\langle f_1 , f_2 \rangle = \int\limits_\R f_1 \bar{f}_2 \, d\lambda ,$$
	where $\lambda$ is the Lebesgue measure on $\R$, and we generalize the notation when it makes sense. We sometimes omit $\R$, when the domain of integration is clear.
	
	For an operator $T\,:\,X\to Y$ between two Banach spaces, $\|T\|_{X\to Y}$ stands for the operator norm, and, sometimes, we omit the subscript simply \text{writing}~$\|T\|$.
	
	Let $\sinc_a\,:\,\C\to\C$ denote the following sine cardinal type function:
	\begin{equation}\label{sinca}
		\sinc_a(z)=\frac{\sin(2\pi a z)}{\pi z}, \qquad z\in\C.
	\end{equation}
	
	\medskip
	
	\subsection{Earlier results}
	We use notation $\mathbb{T} = \{z \in \mathbb{C} : |z| = 1\}$ for the unit circle. Let $m$ denote the Lebesgue measure on $\T$ normalized by $m(\T)=1$. Define the Fourier coefficients of $f\in L^1(\T)$ by 
	$$\widehat{f}(n) = \int\limits_\T f(z) \bar{z}^n \,dm(z),\qquad n\in\Z.$$
	We recall that for $1\leqslant p <+\infty$, a function $f$ on $\T$ is said to belong to the Hardy space $H^p(\D)$ in the unit disk $\D=\{z\in\C:|z|<1\}$ if $f \in L^p(\T)$ and $\widehat{f}(n) = 0$ for all integer $n <0$. We often omit $\D$ and write simply $H^p$. The space $H^p$ is a closed subspace of $L^p(\T)$. Denote by $\mathrm{P}_+$ the orthogonal projection in $L^2(\T)$ to the subspace $H^2$. The classical Toeplitz operator $\mathop{T_\varphi}\,\colon\,H^2\to H^2$ with symbol $\varphi\in L^\infty(\T)$ is defined by 
	$$T_\varphi\; \colon\; f  \mapsto \mathrm{P}_+[\varphi \cdot f],\qquad f \in H^2.$$
	In 1964, A. Brown and P. Halmos \citep{halmos} described basic algebraic properties of Toeplitz operators on $H^2$. In particular, they proved that the Toeplitz operator $T_\varphi$ on $H^2$ with a bounded symbol $\varphi$ satisfies
	$$\left \|T_\varphi \right \|_{H^2\to H^2} = \left \| \varphi \right \|_{L^\infty(\T)},$$
	see Corollary to Theorem 5 in \citep{halmos}. This formula implies that the symbol of a Toeplitz operator on $H^2$ is unique.  
	
	\medskip
	
	For Toeplitz operators on the Paley-Wiener space $\mathrm{PW}_a^2$, the classical treatment of their properties is due to R. Rochberg \citep{RR}. In 1987, he considered boundedness and compactness, as well as Schatten classes $\mathcal{S}^p$ membership. As we mentioned above, he proved that every bounded Toeplitz operator on $\mathrm{PW}_a^2$ admits a bounded symbol. In this paper, we apply his methods to prove a similar result for Toeplitz operators on $\pw$.
	
	\medskip
	
	Toeplitz operators on the Paley-Wiener space are in fact examples of the general truncated Toeplitz operators defined below. A function $\theta \in H^2$ is called an inner function if $|\theta|=1$ $m$-almost everywhere on the unit circle $\T$. With each non-constant inner function $\theta$ we associate the subspace $K_{\theta}^2(\D)=K_\theta^2 = H^2 \ominus \theta H^2$ of $L^2(\mathbb{T})$. Such subspaces are called model subspaces in the unit disk, \citep{nikk}. Denote by $\mathrm{P}_\theta$ the orthogonal projector from $L^2(\T)$ onto $K_{\theta}^2$. A truncated Toeplitz operator $T_\varphi\,\colon\,K_{\theta}^2\to K_{\theta}^2$ with symbol $\varphi \in L^2(\T)$ is densely defined by the following expression
	$$T_\varphi\; \colon\; f  \mapsto \mathrm{P}_\theta[\varphi \cdot f],\qquad f\in K_{\theta}^2 \cap L^\infty(\T).$$
	Toeplitz operators on the Paley-Wiener space are closely related to truncated Toeplitz operators on the model subspace $K^2_{\theta_a}(\C_+)$ of the Hardy space $H^2_+$ in the upper-half plane $\C_+ = \{z\in\C\mid \im z>0\}$ associated with the inner function $\theta_a = e^{2\pi i a z}, a>0$. In fact, $\mathrm{PW}_a^2=\bar{\theta}_a K^2_{\theta_a^2}$, see \citep{nikk}. 
	
	General theory of truncated Toeplitz operators has been pioneered by D. Sarason's paper \citep{saras}, 2007. It plays the same role for truncated Toeplitz operators as the paper of A. Brown and P. Halmos \citep{halmos} does for classical Toeplitz operators. D. Sarason posed several open questions on truncated Toeplitz operators including the problem of the existence of a bounded symbol for a general bounded truncated Toeplitz operator.
	
	\medskip
	
	In 2010, A. Baranov, I. Chalendar, E. Fricain, J. Mashreghi, and D. Timotin \citep{dan}  constructed an inner function $\theta$ and a bounded truncated Toeplitz operator on $K_\theta^2$ that admits no bounded symbol. In 2011, A. Baranov, R. Bessonov, and V. Kapustin \citep{bessonov} characterized inner functions $\theta$ such that every bounded Toeplitz operator on $K_\theta^2$ admits a bounded symbol. In particular, this is the case for so-called one-component inner functions. An inner function $\theta$ is called one-component if the set $\{z\colon |\theta| < \varepsilon\}$ is a connected subset (of the unit disk or the upper half plane of the complex plane) for some $0<\varepsilon <1$. Since the set $\{z\in\C_+\colon |\theta_a(z)|<\varepsilon\}$ is connected for every $0<\varepsilon<1$, this result generalizes the aforementioned theorem by R. Rochberg. 
	
	\medskip	
	
	In 2011, M. Carlsson \citep{carl} proved an estimate similar to what we consider here in the paper. Instead of Toeplitz operators on $\mathrm{PW}_a^2$ he dealt with Wiener-Hopf operators on $L^2[0,2a]$. Following \citep{carl}, define truncated Wiener-Hopf operator $W_{\varphi}$ on $L^2[0,2a]$ with symbol $\varphi\in\sr$ by 
	$$W_\varphi [f](x) = \int\limits_\R \widehat{\varphi}(y) f(x+y)\,dy, \qquad x\in [0,2a],$$
	where $f$ is extended by zero to $\R\setminus [0,2a]$. One can consider more general symbols $\varphi$ including tempered distributions, for simplicity of presentation we limit ourselves by the case $\varphi \in \sr$. M. Carlsson obtained the following estimate
	$$\frac{1}{3}\cdot\left \|\varphi \right \|_{L^\infty(\R)} \leqslant \left \|W_\varphi \right \|_{L^2[0,2a]\to L^2[0,2a]} \leqslant \left \|\varphi \right \|_{L^\infty(\R)},$$
	see Theorem 1.1 in \citep{carl}. This implies 
	$$\frac{1}{3}\cdot\left \|\varphi \right \|_{L^\infty(\R)} \leqslant \left \|T_\varphi \right \|_{\mathrm{PW}_a^2 \to \mathrm{PW}_a^2} \leqslant \left \|\varphi \right \|_{L^\infty(\R)}.$$
	Thus, in the case $p=2$, one can take $c = 1$ in Theorem~\ref{result} of the present paper. In this paper we generalize this result for spaces $\pw$, $1<p<+\infty$.
	
	An extended discussion on truncated Toeplitz operators can be found in survey $\citep{fric}$ by I. Chalendar, E. Fricain and D. Timotin.
	
	\medskip 
	
	In 2011, A. Baranov, R. Bessonov, and V. Kapustin \cite{bessonov} proved that the existence of a bounded symbol for every truncated Toeplitz operator on $K_\theta^2$ is equivalent to the result that every function $f\in H^1 \cap \theta^2 \overline{z H^1}$ admits a weak factorization.
	\begin{Theorem}[\cite{bessonov}, Theorem 2.4]\label{13ex} \textit{Let $\theta$ be an inner function on $\T$. The following assertions are equivalent:
			\begin{enumerate}
				\item any bounded truncated Toeplitz operator on $K_\theta^2(\D)$ admits a bounded symbol;
				\item for any function $f\in H^1(\D) \cap \theta^2 \overline{z H^1(\D)}$ there exist $x_k, y_k \in K_\theta^2(\D)$ with $$\sum_{k\geqslant 0} \left \|x_k\right \|_{L^2(\T)}\cdot\left \|y_k\right \|_{L^2(\T)}<+\infty \;\text{ such that }\;f = \sum_{k\geqslant 0} x_k y_k.$$
		\end{enumerate}}
	\end{Theorem}
	
	Since in some sense the existence of a bounded symbol for every bounded Toeplitz operator on $\pw$, $1<p<+\infty$, will be proved, Theorem \ref{13ex} above allows us to assume that Theorem \ref{result2} of the paper holds true. 
	
	In 1990, K. Dyakonov \citep[Theorem~3]{dyak} proved the strong factorization theorem for non-negative functions, that is, for any $f\in\mathrm{PW}_{2a}^1$, $f\geqslant 0$, there is $g\in\mathrm{PW}_a^2$ such that $f=|g|^2$. These can be utilized to prove the weak factorization theorem in the case $p=2$, see Subsection 1.3 and Example 7.1 in \cite{bessonov}. We note that this approach cannot be easily generalized for arbitrary $p\in (1,+\infty)$, because choosing an outer function with modulus $f^{1/2}$ on $\R$ as $g\in \mathrm{PW}^{2p}_a$, $1\leqslant p<+\infty$, gives $|g|^2\in\mathrm{PW}_{2a}^p$, and not $\mathrm{PW}_{2a}^1$.
	
	\medskip
	
	\subsection{Plan of the proof}
	We outline the structure of the paper. In Section \ref{2} we show that projector $\P_a$ on $L^p(\R)$ is bounded and admits an integral representation with the kernel $\sinc_a$ (see (\ref{sinca})). Further, we show that every Toeplitz operator on $\pw$ also admits an integral representation with $\sinc_a$ kernel. 
	
	Let us describe the plan of the proof of Theorem \ref{result}. Fix some $\varphi\in\sr$. In Section \ref{23} we define three smooth and compactly supported functions, which we use to construct left, central, and right parts of the symbol $\varphi$. Next, given a Toeplitz operator $T_\varphi$ on $\pw$ we construct Toeplitz operators $T_\mathfrak{L}$, $T_\mathfrak{C}$, and $T_\mathfrak{R}$. In Proposition \ref{lcr}, we prove the existence of a universal constant $c>0$ such that
	$$\left \| T_\mathfrak{L} \right \|+\left \| T_\mathfrak{C} \right \|+\left \| T_\mathfrak{R} \right \| \leqslant c\cdot \left \| T_\varphi \right \|.$$
	In the beginning of Section \ref{3}, we start with some preliminaries and prove auxiliary statements. Then we prove the upper bound for the norm of the central part of the symbol. In addition, in Section \ref{4} we define Hankel operators with bounded symbols on the Hardy space in the upper half-plane $H^p_+$ and sketch a proof of the Nehari theorem. Furthermore, we show that any Hankel operator with symbol $\bar{\theta}_a^2 \varphi_*$ such that $\varphi_*\in\sr$ and $\supp \widehat{\varphi}_* \subset \R_+$ corresponds to a Toeplitz operator on $\pw$. Finally, in Section \ref{5} we prove the first result of the present paper.
	
	To conclude, we discuss the plan of the proof for the second result, namely Theorem \ref{result2}, which, in fact, is of its own significance and, in particular, is an example of an implementation of the first result. In Section \ref{6}, we obtain the Toeplitz commutator theorem, that is, if for a special function $\omega$, a bounded operator $T$ on $K_\theta^p$ obeys $\langle T[f],g \rangle=\langle T[\omega f],\omega g \rangle$, then $T$ is a Toeplitz operator on $K_\theta^p$. Finally, in Section \ref{7}, we start with the key statement that the series from Theorem \ref{result2} form a predual space to $\mathcal{T}^p(a)$. Then, the first result and the key statement together yield the second result of the paper, the weak factorization theorem.
	
	\medskip
	
	\section{Toeplitz operators as integral operators. Splitting a symbol}\label{2}
	\subsection{Riesz projector and related operators}
	Let $1 \leqslant p <+\infty$. The Hardy space $H^p_+$ in the upper half-plane $\C_+$ can be defined by
	$$H^p_+ = \cl_{L^p(\R)} \{f \in\sr\mid\supp \widehat{f}\subset \R_+\}.$$
	Let also 
	$$H^p_- = \cl_{L^p(\R)} \{f \in\sr\mid\supp \widehat{f}\subset \R_-\}.$$
	Basic theory of Hardy spaces can be found in \cite{cima}, \cite{garnet}, \cite{koosis}, and \cite{levin}. 
	
	\medskip
	
	Define the Riesz projector $\P_+$ to be the Fourier multiplier associated to symbol $\chi_{\R_+}$, where $\chi_{\R_+}$ is the indicator function of $\R_+$. For $1<p<+\infty$, $\P_+$ extends from $\sr$ to a linear bounded operator on $L^p(\R)$, see e.g., Lecture 19.2 and 19.3 in \citep{levin}. Since $\chi_{\R_+}^2 = \chi_{\R_+}$, we have $\P_+^2 = \P_+$, that is, $\P_+$ operator is a linear bounded projector to $H^p_+$ in $L^p(\R)$. Set 
	$$A_p = \left \| \P_+ \right \|_{L^p(\R)\to L^p(\R)}.$$
	It is known that  
	$$ \begin{array}{cc}
		A_p \leqslant \frac{A}{p-1}, & p \xrightarrow{} 1,\\
		A_p \leqslant Ap, & p \xrightarrow{} +\infty,
	\end{array}$$
	for a universal constant $A > 0$, see \citep{garnet}. Consider an inner function $\theta_a(z) = e^{2\pi i a z}$, $a>0$, $|\theta_a|=1$ almost everywhere on $\R$. Let $\mathcal{U}_t$ be the translation operator $\mathcal{U}_t\,:\,f\mapsto f(\cdot + t)$ and recall that $\P_a$ is the projector to $\pw$.
	
	\begin{Lemma}\label{pa}\textit{
		We have $\left \| \mathbb{P}_a\right \|_{L^p(\R)\to L^p(\R)} \leqslant 2 A_p$.}
	\end{Lemma}
	\begin{proof} 
		Notice that
		$$\chi_{[-a,a]}=\chi_{[-2a,+\infty]}-\chi_{[a,+\infty]}=\mathcal{U}_{2a}[\chi_{\R_+}] - \mathcal{U}_{-a}[\chi_{\R_+}].$$
		By the definition of Fourier transform, $\F^{-1}\mathcal{U}_{a}=\bar{\theta}_a \F^{-1}$ for every $a>0$. Hence,
		\begin{eqnarray*}
			\P_a \;=\; \F^{-1}\chi_{[-a,a]}\F& = & \F^{-1} \chi_{[-2a,+\infty]}\F -  \F^{-1}\chi_{[a,+\infty]}\F
			\\
			&= & \F^{-1} \mathcal{U}_{2a}\chi_{\R_+}\mathcal{U}_{-2a}\F  -  \F^{-1}\mathcal{U}_{-a}\chi_{\R_+}\mathcal{U}_{a}\F
			\\
			&=& \bar{\theta}_a^2\P_+ \theta_a^2  - \theta_a \P_+ \bar{\theta}_a .
		\end{eqnarray*}
		The result follows. 
	\end{proof}
	
	Let $C_0^\infty(\R)$ be the space of all complex-valued smooth functions on $\R$ with compact support. Note that $\P_a(L^p(\R))= \pw$. Indeed, since $\pw$ is a closed subspace of $L^p(\R)$, it is enough to prove that $\P_a(E)\subset \pw$ for some subset $E\subset L^p(\R)$ such that $\cl_{L^p(\R)}E=L^p(\R)$ and $\sar \subset E$. This holds for 
	$$E = \{f\mid \exists g\in C_0^\infty(\R):\; f=\widecheck{g}\;\text{ and }\;a,-a \notin \supp g \}.$$ 
	Whence, the operator $\mathbb{P}_a$ is indeed a bounded projector onto Paley-Wiener space $\pw$. 
	
	\medskip
	
	We now derive an integral formula for $\P_a$. 
	
	\begin{Proposition}\label{proj_pw}\textit{
		For $1<p<+\infty$, the projector $\mathbb{P}_a$ admits the following integral representation:
		\begin{equation}\label{proja}
			\mathbb{P}_a[f](x) = \int\limits_\R \sinc_a(x-y)f(y)\,dy, \qquad f\in L^p(\R).
		\end{equation}}
	\end{Proposition}
	\begin{proof}
		Let us first show that function $\sinc_a\in L^p(\R)$ for every $1< p < +\infty$ (the definition of $\sinc_a$ see in (\ref{sinca})). Indeed, this follows from the estimate
		$$|\sinc_a(x)|\leqslant \frac{1}{\pi  |x|},\qquad x\in\R,$$
		and boundedness of $\sinc_a$ near the origin. Therefore, the integral in \eqref{proja} converges and defines the function on $\R$. Since 
		\begin{equation}\label{eqker}
			\widecheck{\chi}_{[-a,a]}(x) =  \int\limits_{-a}^a e^{2\pi i \xi x} \, d\xi = \frac{e^{2\pi i x a} - e^{-2\pi i x a}}{2\pi i x} = \frac{\sin(2\pi x a)}{\pi x} = \sinc_a(x),
		\end{equation}
		formula \eqref{proja} holds for every $f\in \sr$ by the definition of $\P_a$. Take an arbitrary function $f \in L^p(\R)$ and consider a sequence $\{f_n\}_{n\in\mathbb{N}}\subset\sr$ such that $f_n \to f$ in $L^p(\R)$ as $n\to +\infty$. Then $\P_a[f_n] \to \P_a[f]$ in $L^p(\R)$ and one can choose a subsequence $\{f_{n_k}\}$ such that $\P_a[f_{n_k}](x) \to \P_a[f](x)$ as $n\to +\infty$ for almost every $x\in\R$. On the other hand, 
		$$\P_a[f_{n_k}](x) = \int\limits_\R \sinc_a(x-y)f_{n_k}(y)\,dy	$$
		converges to $\int_\R \sinc_a(x-y)f(y)\,dy$ for every $x\in\R$, by H{\"o}lder's inequality. Hence, \eqref{proja} holds for every $f \in L^p(\R)$.
	\end{proof}
	
	\medskip

	From Proposition \ref{proj_pw} it follows that every Toeplitz operator on $\pw$ with symbol $\varphi\in \mathcal{P}(\R)$ admits the following representation
	\begin{equation}\label{tphi}
		T_\varphi[f](x) = \int\limits_\R \sinc_a(x-y)f(y)\varphi(y)\, dy , \qquad f\in\pw .
	\end{equation}
	Given a function $h$ on $\R$, let $h|_{A}$ denote the restriction of $h$ to a subset $A\subset\R$.
	
	\begin{Lemma}\label{symba}\textit{
		If $\varphi \in \sr$ is such that $\widehat{\varphi}|_{[-2a,2a]}=0$, then $T_\varphi = 0$.}
	\end{Lemma}
	\begin{proof}
		Take $f \in \sar$. By definition, 
		$$\F T_\varphi [f](x) = \chi_{[-a,a]}(x) \int\limits_{\R\setminus [-2a,2a]} \widehat{f}(x-y) \widehat{\varphi}(y) \,dy.$$
		For $x \in [-a,a]$ and $y$ such that $|y|>2a$, we have $|x-y|>a$. Hence, for such $x$ and $y$ we have $\widehat{f}(x-y)=0$ because $\supp \widehat{f} \subset [-a,a]$. 
	\end{proof}
	
	\medskip
	
	\subsection{Splitting procedure. Norm estimates}\label{23}	
	Consider a function $\psi_\mathfrak{L}\in C_0^\infty(\R)$ such that $0\leqslant\psi_\mathfrak{L}\leqslant 1$, $\supp\psi_\mathfrak{L}=[-4,-\tfrac{1}{4}]$ and $\psi_\mathfrak{L}|_{[-2,-\tfrac{1}{2}]}=1$. Set $\psi_\mathfrak{R}(x) = \psi_\mathfrak{L}(-x)$ and define $\psi_\mathfrak{C} = \chi_{[-\frac{1}{2}, \frac{1}{2}]} (1-\psi_\mathfrak{L}-\psi_\mathfrak{R})$. Then $\psi_\mathfrak{L}$, $\psi_\mathfrak{C}$, $\psi_\mathfrak{R}$ are smooth compactly supported functions such that $\psi_\mathfrak{L}+ \psi_\mathfrak{C}+\psi_\mathfrak{R}=1$ on $[-2,2]$, see Figure \ref{fig1} below. 
	
	\begin{figure}[h]
		\centering
		\includegraphics[scale=1, trim={1cm 0 1cm 0},clip]{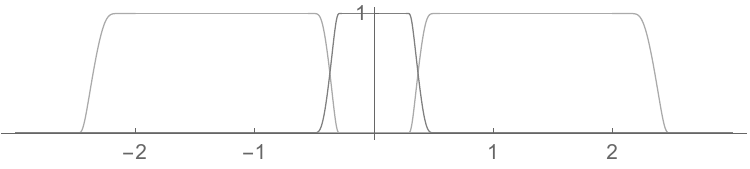}
		\caption{Graphs of functions $\psi_\mathfrak{L}, \psi_\mathfrak{C},\psi_\mathfrak{R}$.}\label{fig1}\ \vspace{-0.4cm}
	\end{figure}
	
	For $a>0$ define $\psi_{\mathfrak{C}, a}\colon x\mapsto \psi_\mathfrak{C}(x/a)$ and $\psi_{\mathfrak{L}, a}$, $\psi_{\mathfrak{R}, a}$ similarly.	Consider a Toeplitz operator $T_\varphi \,\colon\,\pw\to\pw$ with symbol $\varphi\in\sr$. Define $\varphi_\mathfrak{C} = \F^{-1}\psi_{\mathfrak{C},a}\F [\varphi]$ and let $T_\mathfrak{C} = T_{\varphi_\mathfrak{C}}$. Analogously, define $\varphi_\mathfrak{L}$, $\varphi_\mathfrak{R}$, $T_{\mathfrak{L}}$, $T_\mathfrak{R}$ using functions $\psi_\mathfrak{L}$, $\psi_\mathfrak{R}$. We call $T_\mathfrak{L}$, $T_{\mathfrak{C}}$, $T_\mathfrak{R}$ the left, central, and right parts of $T_\varphi$, respectively. 

	\begin{Proposition}\label{lcr}\textit{
		Let $1<p<+\infty$. Consider a Toeplitz operator $T_\varphi$ on $\pw$ with symbol $\varphi \in \sr$. We have $T_{\varphi} = T_\mathfrak{L}+T_\mathfrak{C}+T_\mathfrak{R}$ and
		$$ c\Bigl(\left \|T_\mathfrak{L}\right \| + \left \|T_\mathfrak{C}\right \| + \left \|T_\mathfrak{R}\right \|\Bigr) \leqslant \left \| T_\varphi \right \| \leqslant \left \|T_\mathfrak{L}\right \| + \left \|T_\mathfrak{C}\right \| + \left \|T_\mathfrak{R}\right \|,$$
		for a universal constant $c > 0$.}
	\end{Proposition}
	\begin{proof}
		Since $\psi_{\mathfrak{L},a}+\psi_{\mathfrak{C},a}+\psi_{\mathfrak{R},a}=1$ on $[-2a,2a]$, we have $\widehat{\varphi}=\widehat{\varphi}_\mathfrak{L}+\widehat{\varphi}_\mathfrak{C}+\widehat{\varphi}_\mathfrak{R}$ on $[-2a,2a]$. Hence, we have $T_{\varphi} = T_\mathfrak{L}+T_\mathfrak{C}+T_\mathfrak{R}$ by Lemma \ref{symba}, and thus
		$$ \left \| T_\varphi \right \|\leqslant\left \|T_\mathfrak{L}\right \| + \left \|T_\mathfrak{C}\right \| + \left \|T_\mathfrak{R}\right \| .$$
		Let us check the opposite inequality. Take $f \in \sar$. By the Fubini–Tonelli Theorem and \eqref{tphi}, we have
		\begin{eqnarray*}
			T_\mathfrak{C}[f](x)& = & \int\limits \sinc_a(x-y)f(y)\cdot \left[ \int\limits \varphi(y-t)\widecheck{\psi}_{\mathfrak{C},a}(t)  \,dt \right]    \,dy 
			\\
			&=& \int\limits \widecheck{\psi}_{\mathfrak{C},a}(t)  \int\limits \sinc_a(x-y) f(y) \varphi(y-t)   \,dy    \,dt    
			\\
			&= & \int\limits \widecheck{\psi}_{\mathfrak{C},a}(t) \int\limits \sinc_a(x-t-\xi)f(\xi+t)\varphi(\xi) \,d\xi    \,dt 
			\\
			&= &\int\limits \widecheck{\psi}_{\mathfrak{C},a}(t) \cdot \mathcal{U}_{-t}T_\varphi \mathcal{U}_t[f](x)    \,dt .
		\end{eqnarray*}
		Set $p(x)= |\widecheck{\psi}_{\mathfrak{C},a}(x)| /\left \|\widecheck{\psi}_{\mathfrak{C},a} \right \|_{L^1(\R)}$, then $\int_\R p(x)\,dx = 1$. Jensen's inequality gives 
		$$\Phi\Biggl(\,\int\limits_\R h(x)p(x)\,dx \Biggr) \leqslant \int\limits_\R \Phi(h(x)) p(x)\,dx, $$
		for every convex function $\Phi\, \colon\, \R \to \R_+$ and every $h$ such that $h p\in L^1(\R)$. Choosing $\Phi = |x|^p$, we obtain
		\begin{eqnarray*}
			\left \| T_\mathfrak{C}[f] \right \|_{L^p(\R)}^p &=& \int\limits \left| \int\limits \widecheck{\psi}_{\mathfrak{C},a}(t) \cdot \mathcal{U}_{-t}T_\varphi \mathcal{U}_t[f](x)    \,dt   \right|^p \,dx 
			\\
			&= & \left \| \widecheck{\psi}_{\mathfrak{C},a} \right \|_{L^1(\R)}^p  \int\limits \left| \int\limits  |\mathcal{U}_{-t}T_\varphi \mathcal{U}_t[f](x) |\cdot p(t)  \,dt   \right|^p \,dx
			\\
			&\leqslant & \left \| \widecheck{\psi}_{\mathfrak{C},a} \right \|_{L^1(\R)}^{p} \int\limits  \int\limits   |\mathcal{U}_{-t}T_\varphi \mathcal{U}_t[f](x) |^p \cdot p(t) \,dt    \,dx 
			\\
			&= & \left \| \widecheck{\psi}_{\mathfrak{C},a} \right \|_{L^1(\R)}^{p-1} \int\limits |\widecheck{\psi}_{\mathfrak{C},a}(t)| \int\limits  |\mathcal{U}_{-t}T_\varphi \mathcal{U}_t[f](x) |^p  \,dx    \,dt
			\\
			&=& \left \| \widecheck{\psi}_{\mathfrak{C},a} \right \|_{L^1(\R)}^{p-1} \int\limits |\widecheck{\psi}_{\mathfrak{C},a}(t)|\cdot  \left \|\mathcal{U}_{-t}T_\varphi \mathcal{U}_t [f] \right \|_{L^p(\R)}^p  \,dt 
			\\
			&\leqslant&  \left \| T_\varphi \right \|^p \cdot \left \| f \right \|_{L^p(\R)}^p \cdot \left \| \widecheck{\psi}_{\mathfrak{C},a} \right \|_{L^1(\R)}^p  .
		\end{eqnarray*}
		Hence, $\left \| T_\mathfrak{C} \right \| \leqslant \left \| T_\varphi \right \|\cdot\left \| \widecheck{\psi}_{\mathfrak{C},a} \right \|_{L^1(\R)}$.
		Similar arguments apply to $T_\mathfrak{L}$, $T_\mathfrak{R}$ and give us the estimate
		$$\left \|T_\mathfrak{L}\right \| + \left \|T_\mathfrak{C}\right \| + \left \|T_\mathfrak{R}\right \| \leqslant \left(\left \| \widecheck{\psi}_{\mathfrak{L},a} \right \|_{L^1(\R)}+\left \| \widecheck{\psi}_{\mathfrak{C},a} \right \|_{L^1(\R)}+\left \| \widecheck{\psi}_{\mathfrak{R},a} \right \|_{L^1(\R)}\right) \cdot \left \| T_\varphi \right \| .$$
		Observe that the constant in the right hand side does not depend on $a$ because
		$$\left \| \widecheck{\psi}_{\mathfrak{C},a} \right \|_{L^1(\R)} = \left \| \widecheck{\psi}_{\mathfrak{C}} \right \|_{L^1(\R)} $$
		and similar identities hold for $\widecheck{\psi}_{\mathfrak{L},a}$, $\widecheck{\psi}_{\mathfrak{R},a}$.
	\end{proof}
	
	\medskip
	
	\section{Reproducing kernels. Central part of a symbol}\label{3}	
	In Section \ref{2} we prove that $\P_a(L^p(\R)) = \pw$ for $1<p<+\infty$. In addition, Proposition \ref{proj_pw} says that for every function $f\in\pw$ we have
	\begin{equation}\label{id}
		f(x)=\P_a[f](x)=\int\limits_\R \sinc_a(x-y)f(y)\,dy,
	\end{equation}
	almost everywhere on $\R$. Note that the right hand side is an entire function with respect to $x$. This follows from the result that the integral $$\int\limits_\R \frac{\partial}{\partial z} \sinc_a(z-y)f(y)\,dy=\int\limits_\R \frac{2 a \cos (2\pi a (z-y))- \sinc_a(z-y)}{z-y}f(y)\,dy$$
	converges uniformly in a neighborhood of any point $z\in\C$. This shows that any function $f\in\pw$ can be naturally identified with an entire function using \eqref{id}. In other words, for every $f\in\pw$ one can find an entire function $g\,\colon\,\C\to\C$ such that $g\in L^p(\R)$ and $f = g$ almost everywhere on $\R$. In particular, for every $z\in\C$ and $f\in\pw$ the value $f(z)$ is well defined.
	
	\begin{Lemma}\label{value_point}\textit{
		Let $1<p<+\infty$ and $\frac{1}{p}+\frac{1}{q}=1$. For each $z \in \C$ the linear functional $\phi_z\,\colon\, f \mapsto f(z)$ on $\pw$ is bounded and 
		$$\phi_z(f)=\int\limits_\R \sinc_a(z-y) f(y)\,dy,\qquad f \in \mathrm{PW}_a^p. $$
		Moreover, for $x\in\R$ we have $\left\| \phi_x\right \|\leqslant  \left \| \sinc_a \right \|_{L^q(\R)}$.}
	\end{Lemma}
	\begin{proof}
		By definition, we have
		$$\phi_z(f) = f(z) = \int\limits_\R \sinc_a(z-y) f(y)\,dy, \qquad f\in\pw .$$
		Then, by H{\"o}lder's inequality for every $f\in\pw$ we have 
		$$|\phi_z(f)| \leqslant \left \| \mathcal{U}_{-z}[\sinc_a] \cdot f \right \|_{L^1(\R)} \leqslant \left \| \mathcal{U}_{-z}[\sinc_a] \right \|_{L^q(\R)}\cdot\left \|f \right \|_{L^p(\R)}.$$
		It follows that $\phi_z$ is bounded and $\left\| \phi_z\right \|\leqslant  \left \| \mathcal{U}_{-z}[\sinc_a] \right \|_{L^q(\R)}$. In particular, if $x\in\R$, then $\left\| \phi_x\right \|\leqslant  \left \| \sinc_a \right \|_{L^q(\R)}$.
	\end{proof}
	
	Now, we obtain an upper bound for the norm of the central part of a symbol.
	
	\begin{Proposition}\label{boundsymb2}\textit{
		Let $1<p<+\infty$. Consider a Toeplitz operator $T_{\varphi}$ on $\pw$ with symbol $\varphi \in\sr$. Let $T_\mathfrak{C}$ be its central part constructed in Section \ref{23}. Then we have
		$$\left \| \varphi_\mathfrak{C} \right \|_{L^\infty(\R)} \leqslant c_{p}\cdot \left \| T_\mathfrak{C} \right \|_{\pw\to\pw},$$
		for some constant $c_{p}>0$ depending only on $p$.}
	\end{Proposition}
	\begin{proof}
		Take $\varepsilon= \tfrac{a}{8}$ and fix some $x\in\R$. From formula \eqref{eqker} we see that $\supp\F[\sinc_\varepsilon(\cdot)]\subset [-\varepsilon, \varepsilon]$, therefore, $\sinc_{\varepsilon}\in\pw$. Recall that $\supp \widehat{\varphi}_{\mathfrak{C}} = [-\tfrac{a}{2},\tfrac{a}{2}]$, hence the support of 
		$$\F[\varphi_\mathfrak{C}\cdot\mathcal{U}_{-x}[\sinc_{\varepsilon}]]=(\psi_{\mathfrak{C},a}\widehat{\varphi}) * (\chi_{[-\varepsilon,\varepsilon]}e^{-2\pi i x \xi})$$
		is in $[-a,a]$ by properties of convolution ($\supp f * g \subset \supp f + \supp g$). We have
		\begin{eqnarray*}
			\phi_x(T_{\mathfrak{C}}\mathcal{U}_{-x}[\sinc_{\varepsilon}])& = & T_{\mathfrak{C}}\mathcal{U}_{-x}[\sinc_{\varepsilon}](x)
			\\
			&=&\P_a[\varphi_\mathfrak{C}\cdot\mathcal{U}_{-x}[\sinc_{\varepsilon}]](x)
			\\
			&=&\varphi_\mathfrak{C}(x)\cdot\mathcal{U}_{-x}[\sinc_{\varepsilon}](x)
			\\
			&= &\varphi_\mathfrak{C}(x)\cdot\sinc_{\varepsilon}(0)
			\\
			&=&2\varepsilon\cdot \varphi_\mathfrak{C}(x).
		\end{eqnarray*}
		By Lemma \ref{value_point}, we have $\phi_x \in (\mathrm{PW}_a^p)^*$, therefore
		\begin{eqnarray*}
			|\varphi_\mathfrak{C}(x)|& \leqslant &\frac{1}{2\varepsilon} \left\|\phi_x \right \| \cdot \left\|T_{\mathfrak{C}}\mathcal{U}_{-x}[\sinc_\varepsilon] \right \|_{L^p(\R)}
			\\
			&\leqslant &\frac{1}{2\varepsilon} \left\|\sinc_a \right \|_{L^q(\R)} \cdot\left \| T_{\mathfrak{C}}\right \|\cdot \left\|\sinc_\varepsilon \right \|_{L^p(\R)}.
		\end{eqnarray*}
		Observe that the latter product of norms does not depend on $a$:
		$$\frac{1}{2\varepsilon} \left\|\sinc_a \right \|_{L^q(\R)} \cdot\left\|\sinc_\varepsilon \right \|_{L^p(\R)} = 4 \left\|\sinc_1 \right \|_{L^q(\R)} \cdot\left\|\sinc_{1/8} \right \|_{L^p(\R)}.$$
	\end{proof}
	
	\begin{Lemma}\label{cub}\textit{
		Let $1<p<+\infty$ and $\tfrac{1}{p}+\tfrac{1}{q}=1$. For a universal constant $c>0$,
		$$\left\|\sinc_1 \right \|_{L^q(\R)} \cdot\left\|\sinc_{1/8} \right \|_{L^p(\R)} \leqslant c\cdot \biggl(p+ \frac{1}{p-1}\biggl).$$
		}
	\end{Lemma}
	\begin{proof}
		We have
		$$\left\|\sinc_{1/8} \right \|_{L^p(\R)} = 8^{-\frac{1}{q}} \left\|\sinc_1 \right \|_{L^p(\R)}\leqslant \left\|\sinc_1 \right \|_{L^p(\R)}.$$
		Clearly, $|\sinc_1(x)| \leqslant 2$ for $|x|\leqslant\frac{1}{2\pi}$ and $|\sinc_1(x)| \leqslant \frac{1}{\pi |x|}$ for $|x|>\frac{1}{2\pi}$. Then, we obtain
		$$\left\|\sinc_1 \right \|_{L^q(\R)}^q\leqslant \frac{2^q}{\pi}  + \frac{2}{\pi}\int\limits_{1/2}^{+\infty} \frac{dx}{x^q} = \frac{2^q}{\pi}\biggl(1+ \frac{1}{q-1}\biggl). $$
		Then,
		$$\Biggl(\frac{2^q}{\pi}\biggl(1+ \frac{1}{q-1}\biggl)\Biggr)^{\frac{1}{q}} \cdot \Biggl(\frac{2^p}{\pi}\biggl(1+ \frac{1}{p-1}\biggl)\Biggr)^{\frac{1}{p}} = \frac{4}{\pi} \biggl(1+ \frac{1}{q-1}\biggl)^{\frac{1}{q}} \cdot \biggl(1+ \frac{1}{p-1}\biggl)^{\frac{1}{p}}, $$
		and, by Bernoulli's inequality, we get
		\begin{align*}
			\biggl(1+ \frac{1}{q-1}\biggl)^{\frac{1}{q}} \cdot \biggl(1+ \frac{1}{p-1}\biggl)^{\frac{1}{p}} &\leqslant  \biggl(1 + \frac{1}{(q-1)q}\biggl) \cdot \biggl(1 + \frac{1}{(p-1)p}\biggl) 
			\\
			&\leqslant \biggl(1 + \frac{1}{q-1 }\biggl) \cdot \biggl(1 + \frac{1}{(p-1)p}\biggl)
			\\
			& =  p \cdot \biggl(1 + \frac{1}{(p-1)p}\biggl)
			\\
			&= p+ \frac{1}{p-1}.
		\end{align*}\nobreak{To sum up, one can pick $c =\tfrac{4}{\pi}$.}
	\end{proof}
	
	\medskip
	
	\section{Nehari Theorem. Right and left parts of a symbol}\label{4}
	\subsection{Hankel operators on the Hardy space. Nehari Theorem}\label{4.1}
	A Hankel operator $H_\varphi \,\colon\, H^2 \to \overline{z H^2}$ with symbol $\varphi \in L^2(\T)$ can be densely defined by 
	$$H_\varphi\; \colon\; f  \mapsto \mathrm{P}_-[\varphi \cdot f], \qquad f\in H^2\cap L^\infty(\T),$$
	where $\mathrm{P}_- = I - \mathrm{P}_+$. Consider $p$ such that $1<p<+\infty$. Similarly, one can define Hankel operator $H_\varphi \,\colon\,H^p_+\to H^p_-$ with symbol $\varphi \in L^\infty(\R)$ by 
	$$H_\varphi\; \colon\; f  \mapsto \P_-[\varphi \cdot f], \qquad f\in H^p_+,$$
	where $\P_- = I - \P_+$, $I$ being the identity operator on $L^p(\R)$. For an introduction to the theory of Hankel operators, see the monograph \cite{peller} by V. Peller. The following theorem, which characterizes bounded Hankel operators on $H^2$, is due to Z. Nehari.
	\begin{Theorem}[\cite{peller}, Theorem 1.3] \textit{ Let $\varphi \in L^2(\T)$. The following statements are equivalent:
		\begin{enumerate}
			\item $H_\varphi$ is bounded on $H^2$;
			\item there exists $\psi\in L^\infty(\T)$ such that $H_\psi = H_\varphi$ and $\left \| \psi \right \|_{L^\infty(\T)} = \left \| H_\varphi \right \|_{H^2 \to \overline{z H^2}}$.
		\end{enumerate}}
	\end{Theorem}
	
	The following theorem can be proved in the same way as Nehari's theorem.
	
	\begin{Theorem}\label{nehari} \textit{Let  $1<p<+\infty$ and let $\varphi \in L^\infty(\R)$. Then there exists a function $\psi\in L^\infty(\R)$ such that $H_\psi=H_\varphi$ and, moreover, $\left \| \psi \right \|_{L^\infty(\R)} \leqslant \left \| H_\varphi \right \|_{H^p_+\to H^p_-}$.}
	\end{Theorem}
	
	We give a sketch of the proof of this result in \hyperref[appendix]{Appendix}.
	
	\medskip
	
	\subsection{Analytic Toeplitz operators on $\mathrm{PW}_a^p$ as Hankel operators}\label{4.2}
	We say that a Toeplitz operator $T_\varphi$ with symbol $\varphi \in\sr$ is called an analytic operator if $\supp \widehat{\varphi}\subset \R_+$. One can easily check that that for every $1<p<+\infty$ and for every $a>0$ we have
	$$\P_a=\theta_a{\P}_-\bar{\theta}_a^2{\P}_+ \theta_a.$$
	This formula will be used in the proof of Lemma \ref{ht} below.
	\begin{Lemma}\label{ht}\textit{
		Let $1<p<+\infty$ and let $\varphi \in \sr$ be such that $\supp \widehat{\varphi} \subset \R_+$. Then
		\begin{equation}\label{www}
			H_{\bar{\theta}_a^2 \varphi}=\bar{\theta}_a T_\varphi \theta_a \P_- \bar{\theta}_a^2 .
		\end{equation}}
	\end{Lemma}
	\begin{proof}
		Note that for any function $g \in H^p_+$, there are functions $g_1 \in\pw$, $g_2 \in H^p_+$ such that $g = \theta_a g_1 + \theta^2_a g_2$. We have
		$$H_{\bar{\theta}_a^2 \varphi}[g]=\P_-[\bar{\theta}_a \varphi g_1 + \varphi g_2] =H_{\bar{\theta}_a^2 \varphi}[\theta_a  g_1 ],$$
		because $\varphi g_2 \in H^p_+$. We also have 
		$$\bar{\theta}_a T_\varphi \theta_a \P_- \bar{\theta}_a^2 [g] = \bar{\theta}_a T_\varphi \theta_a \P_- [\bar{\theta}_a g_1 + g_2]= \bar{\theta}_a T_\varphi [g_1].$$
		On the other hand, taking into account \eqref{www}, we obtain
		$$\bar{\theta}_a T_\varphi[g_1]= \bar{\theta}_a \P_a[\varphi g_1]
		= \P_{-}\bar{\theta}_a^2\P_{+}[\theta_a \varphi g_1]
		=  \P_{-}[\bar{\theta}_a\varphi g_1]
		= H_{\bar{\theta}_a^2 \varphi} [\theta_a g_1].$$
		This completes the proof.
	\end{proof}
	
	\medskip
	
	\section{Existence of a bounded symbol}\label{5}
	Now, we prove the first result, namely Theorem \ref{result}. Every Toeplitz operator $T_\varphi$ on $\pw$, $1<p<+\infty$, with symbol $\varphi \in \sr$ admits a bounded symbol $\psi$ such that 
	$$\left \| \psi \right \|_{L^\infty(\R)} \leqslant c \biggl(p+\frac{1}{p-1}\biggr) \cdot  \left \| T_\varphi \right \|_{\pw\to\pw},$$
	for a universal constant $c > 0$. 
	
	\begin{proof}
		Define operators $T_\mathfrak{L}$, $T_\mathfrak{C}$, $T_\mathfrak{R}$ as in Section \ref{23}. By Proposition \ref{lcr} we have
		$$ \left \|T_\mathfrak{L}\right \| + \left \|T_\mathfrak{C}\right \| + \left \|T_\mathfrak{R}\right \| \leqslant c\cdot \left \| T_\varphi \right \|,$$
		for a universal constant $c > 0$. By Proposition \ref{boundsymb2} we have
		$$\left \| \varphi_\mathfrak{C} \right \|_{L^\infty(\R)} \leqslant c_{p}\cdot \left \| T_\mathfrak{C} \right \|,$$
		for some constant $c_{p}>0$ depending only on $p$. We now prove an upper bound for the left and right parts of Toeplitz operators. By the Nehari Theorem (see Theorem \ref{nehari}), there exists $\psi_r \in L^\infty(\R)$ such that $H_{\psi_r} = H_{\bar{\theta}_a^2 \varphi_\mathfrak{R}}$, and, moreover,
		$$\left\|\psi_r\right \|_{L^\infty(\R)}\leqslant\left\| H_{\psi_r} \right \|=\left \| \bar{\theta}_a T_\mathfrak{R} \theta_a \P_- \bar{\theta}_a^2 \right \| \leqslant A_p \left \| T_\mathfrak{R} \right \|,$$
		where we used the result that $\left \| \P_- \right \|=\left \| \P_+ \right \|= A_p$. 
		
		Next, we show that $T_\mathfrak{R} = T_{\theta_a^2 \psi_r}$. Since $H_{\psi_r} = H_{\bar{\theta}_a^2 \varphi_\mathfrak{R}}$, we have $H_{\psi_r}[\theta^2_a f] = H_{\bar{\theta}_a^2 \varphi_\mathfrak{R}}[\theta^2_a f]=0$ for every $f \in H^p_+$. Therefore, 
		$$\P_+[\psi_r \theta_a^2 f] = \psi_r \theta_a^2 f - H_{\psi_r}[\theta^2_a f]=\psi_r \theta_a^2 f, \quad f \in H^p_+.$$ 
		Let $h\in\pw$ and let $f = \theta_a h$. Then $f\in H^p_+$ and we have
		\begin{align*}
			T_{\theta_a^2 \psi_r}[h] &= \P_a[\theta_a^2 \psi_r h]= \theta_a \P_- \bar{\theta}_a^2 \P_+ [\theta^2_a \psi_r f] = \\
			&= \theta_a \P_- [\psi_r f] = \theta_a H_{\psi_r} [f] = \theta_a H_{\bar{\theta}_a^2 \varphi_\mathfrak{R}}[f].	
		\end{align*}
		By Lemma \ref{ht}, we have $\theta_a H_{\bar{\theta}_a^2 \varphi_\mathfrak{R}}[f] = T_\mathfrak{R}\theta_a \P_- [\bar{\theta}_a^2 f]=T_\mathfrak{R}\theta_a \P_- [\bar{\theta}_a h] = T_\mathfrak{R}[h]$, so the claim follows.
		
		Similarly, there exists $\psi_l \in L^\infty(\R)$ such that
		$$\left \| \psi_l \right \|_{L^\infty(\R)} \leqslant  A_p \left \| T_\mathfrak{L}\right \|\;\text{ and }\; T_\mathfrak{L} = T_{\bar{\theta}_a^2 \psi_l}.$$
		Setting $\psi =\bar{\theta}_a^2 \psi_l +\varphi_\mathfrak{C} + \theta_a^2 \psi_r$ we obtain 
		$$T_\varphi = T_\mathfrak{L}+T_\mathfrak{C}+T_\mathfrak{R}= T_{\bar{\theta}_a^2 \psi_l}+T_\mathfrak{C}+T_{\theta_a^2 \psi_r} = T_\psi,$$
		by Proposition \ref{lcr}. Since 
		\begin{align*}
			\left \|\psi \right \|_{L^\infty(\R)} &\leqslant \left \|\psi_l \right \|_{L^\infty(\R)}+ \left \|\varphi_\mathfrak{C} \right \|_{L^\infty(\R)}+\left \|\psi_r \right \|_{L^\infty(\R)}  \\
			&\leqslant  A_p \left \| T_\mathfrak{L}\right \|+ c_p \left \| T_\mathfrak{C}\right \|+ A_p \left \| T_\mathfrak{R}\right \| \\
			&\leqslant \tilde{c}\cdot(2A_p +c_p) \left \| T_\varphi \right \|,
		\end{align*}
		we have 
		$$	\left \|\psi \right \|_{L^\infty(\R)} \leqslant c \cdot \biggl(p + \frac{1}{p-1}\biggr) \left \| T_\varphi \right \|,$$
		by Lemma \ref{cub} and the estimate for the Riesz projector norm from Section \ref{2}. The theorem is proved.
	\end{proof}
	
	\medskip
	
	\section{Characteristic property of Toeplitz operators}\label{6}
	We now turn to the proof of the second result. In this section, we show the Toeplitz commutator theorem we mentioned in the abstract. 
	
	Let $1< p < +\infty$ and $q$ be the H{\"o}lder conjugate of $p$. A bounded analytic function $\theta\,:\,\C_+\to\C$ is an inner function in the upper half-plane $\C_+$ if 
	$$\lim\limits_{y\to 0+}|\theta(x+iy)|=1,\quad\text{ a.e. } x\in\R.$$
	Here, a.e. means almost everywhere with respect to the Lebesgue measure on $\R$. The model subspace $K_\theta^p$ is defined by
	$$K_\theta^p = H_+^p\cap \theta H_-^p.$$
	The involution on $K_\theta^p$ is given by $\tilde{f}=\theta\bar{f}$. Notice that from this point on, the notation $K_\theta^p$ corresponds to the model subspace in the upper half-plane $\C_+$, and when it is needed we write $K_\theta^p(\C_+)$ and $K_\theta^p(\D)$ for the model subspace in the upper half-plane $\C_+$ and the model subspace in the unit disk $\D$, respectively. More information about the model subspaces can be found \text{in \cite{dyak}, \cite{nikk}.}
	
	Fix the following function $\sigma(x) = (x+i)^{\frac{2}{p}},\,x\in\R$. Let $\mathbb{P}_\theta$ be a projection operator from $L^p(\R)$ onto $K_\theta^p$. A Toeplitz operator $T_\varphi$ on $K_\theta^p$ with symbol $\varphi \in \sigma (L^p(\R)+L^q(\R))$ is defined by 
	$$T_\varphi\,\colon\, f\mapsto \P_\theta[\varphi\cdot f],\qquad f\in K_\theta^p.$$
	Let $H(\C_+)$ denote the space of holomorphic functions in $\C_+$. The space $H^\infty_+$ is the space of bounded functions from $H(\C_+)$, with the standard supremum norm. A Toeplitz operator can be thought of in two ways: either as an unbounded operator $K_\theta^p\to K_\theta^p$ whose domain contains $K_\theta^p\cap H^\infty_+$, or as an operator $K_\theta^p \to H(\C_+)$, continuous relative to the weak topology of $K_\theta^p$ and the topology of locally uniform convergence of $H(\C_+)$. We study bounded Toeplitz operators on $K_\theta^p$. The space of all bounded Toeplitz operators on $K_\theta^p$ is denoted by $\mathcal{T}^p(\theta)$.
	
	In 2007, D. Sarason \cite{saras} proved the characteristic property of Toeplitz operators.
	\begin{Theorem}[\cite{saras}, Theorem 8.1] \textit{
		A bounded operator $A$ on $K_\theta^2(\mathbb{D})$ is a truncated Toeplitz operator if and only if the following holds:
		$$\left\langle A[x], x \right \rangle = \left\langle A[zx], zx \right \rangle, \qquad x, zx \in K_\theta^2(\mathbb{D}).$$}
	\end{Theorem}
	Recall that for $z\in\C_+$ the following conformal map 
	$$\omega\,\colon\,z\mapsto\frac{z-i}{z+i}\in\mathbb{D}$$
	sends $\C_+$ to the unit disk $\D$ and $\omega(z)|_{z\in\mathbb{R}}\in\mathbb{T}$. Define the operator $\mathbb{U}\,\colon\,L^p(\T)\to L^p(\R)$~by
	\begin{equation}\label{prem1}
		\mathbb{U}[f](x) = \left(\frac{1}{\pi(x+i)^2}\right)^{\frac{1}{p}}\,f\left(\omega(x)\right), \quad x\in\R.
	\end{equation}
	Note that $\mathbb{U}$ is an isometric isomorphism, see Chapter 6 in \citep{NikVol1}. Clearly, $\mathbb{U}[zf] = \omega \, \mathbb{U}[f]$ for any $f\in L^p(\T)$. This identity and the theorem above immediately allows us to formulate the following hypothesis.
	\begin{Corollary}\label{con3} \textit{
		A bounded operator $T\,:\,K_\theta^p\to K_\theta^p$ is a Toeplitz operator on $K_\theta^p$ if and only if the following holds: 
		$$\left\langle T[f], g \right \rangle = \left\langle T[\omega f], \omega g \right \rangle, \qquad f, \omega f \in K_\theta^p,\, g, \omega g \in K_\theta^q.$$}
	\end{Corollary} 
	All remaining subsections of this section are devoted to proving the corollary above.
	
	\medskip
	
	\subsection{Conjugate kernel}
	Let $1< p < +\infty$ and $q$ be the H{\"o}lder conjugate of $p$. The Riesz projector $\mathbb{P}_+ \,\colon L^p(\R)\to H_+^p$ may be written as a Cauchy type integral
	$$\mathbb{P}_+[f](z) = \frac{1}{2\pi i}\int\limits_\R \frac{f(x)}{x-z}\,dx,\qquad z\in\C_+.$$
	Given function $h_z \in H^p_+$ such that
	$$h_z(x) = \frac{1}{2\pi i}\,\frac{1}{\bar{z}-x},\qquad x\in\R, z\in\C_+.$$
	For any $f \in H^p_+$ we have $\left\langle f, h_z \right \rangle = f(z),\; \mathrm{Im}z>0$, see Lecture 19.2 in \cite{levin}. Hence, for any $f \in K_\theta^p$, we get
	$$f(z)=\left\langle f, \mathbb{P}_\theta [h_z] \right\rangle,\qquad \mathrm{Im} z>0.$$
	We define $k_{\theta,z} = \mathbb{P}_\theta [h_z]$, and these $k_{\theta,z}$ is the reproducing kernel for the space $K_\theta^p$. 
	
	Notice that for any $g\in H^q_+$,
	$$\left \langle \mathbb{P}_+[\bar{\theta}h_z], g \right\rangle = \left \langle \bar{\theta}h_z, g \right\rangle = \overline{\left \langle \theta g, h_z \right\rangle} = \overline{\theta(z)}\cdot \overline{g(z)}= \overline{\theta(z)}\left \langle h_z, g \right\rangle,$$
	then $\mathbb{P}_+[\bar{\theta}h_z] = \overline{\theta(z)} h_z$. The projection operator $\mathbb{P}_\theta$ onto $K_\theta^p$ can be rewritten via the Riesz projector, $\mathbb{P}_\theta = \mathbb{P}_+ - \theta \mathbb{P}_+ \bar{\theta}$, whence
	$$k_{\theta,z}(x) = \frac{1}{2\pi i}\,\frac{1-\overline{\theta(z)}\theta(x)}{\bar{z}-x}, \qquad x\in\R, z\in\C_+.$$
	Therefore, the conjugate kernel $\tilde{k}_{\theta,z}$ defined by $\tilde{k}_{\theta,z} = \theta \bar{k}_{\theta, z}$ can be represented in the form
	$$\tilde{k}_{\theta,z}(x) = \frac{1}{2\pi i}\,\frac{\theta(x)-\theta(z)}{x-z}, \qquad x\in\R, z\in\C_+.$$
	
	We remark that $\omega$ is an inner function in the upper half plane $\C_+$. Now, we prove the following auxiliary lemma, and then we show what the condition $f,\omega f\in K_\theta^p$ from Corollary \ref{con3} means in terms of the conjugate kernel $\tilde{k}_{\theta,i}$.
	
	\begin{Lemma}\label{sp} \textit{
		The space $K_\omega^p$ is a linear span of the function $h_i$, in other words 
		$$K_\omega^p = \mathrm{span}\,\left\langle \frac{1}{x+i} \right\rangle .$$}
	\end{Lemma}
	\begin{proof}
		Clearly, $\left\langle h_i, \omega h \right\rangle = \overline{\omega(i) h(i)} = 0$ for any $h \in H^q_+$, hence $\mathrm{span}\langle h_i\rangle \subset H^p_+ \cap \omega {H^p_-}$. Conversely, it remains to prove that if $f \in K^p_\omega$, then there exists a constant $c\in\C$ such that $f = c\cdot h_i$. We show that for some constant $c$, one has 
		$$\omega \overline{\left(f - \frac{c}{x+i}\right)} \in H^p_+\;\text{ and }\;\overline{\omega} \left(f - \frac{c}{x+i}\right) \in H^p_+.$$
		It is easy to see that $\omega \bar{f} \in H^p_+$ and $\omega \bar{h}_i \in H^p_+$, so $\omega\overline{(f-ch_i)}\in H_+^p$ for any $c\in\C$. By the inner-outer factorization theorem (see e.g., Theorem 3.2.4 in \cite{cima} or Corollary 5.7 in \cite{garnet}), we have that there exists an outer function $\mathfrak{F}\in H^p_+$ and an inner function $\mathfrak{I}$ such that
		$$g(x) = f(x) - \frac{2i f(i)}{x+i} = \omega(x)\,\mathfrak{I}(x)\cdot\mathfrak{F}(x),$$
		and $g(i)=0$. Thus, we get
		$$\frac{x+i}{x-i} \left(f - \frac{2i f(i)}{x+i}\right) = \mathfrak{I}\cdot\mathfrak{F}\in H^p_+.$$
	\end{proof}

	\begin{Proposition}\label{l3} \textit{
		Let $f \in K_\theta^p$, then $\omega f \in K_\theta^p$ if and only if $\left\langle  f, \Tilde{k}_{\theta,i}\right\rangle  = 0$.}
	\end{Proposition}
	\begin{proof}
		Observe that $\omega f \in H^p_+$ and for any $x\in\R$ we have
		$$\omega(x) \Tilde{k}_{\theta,i}(x) = \frac{1}{2\pi i}\,\frac{\theta(x)-\theta(i)}{x-\bar{i}} = (\theta(i)-\theta(x)) \cdot h_i(x).$$
		If $\omega f \in K_\theta^p$ we have
		$$\left\langle  f, \Tilde{k}_{\theta,i}\right\rangle = \left\langle \omega f,  \omega \Tilde{k}_{\theta,i}\right\rangle  = \overline{\theta(i)}\left\langle  \omega f,   h_i \right\rangle -  \left\langle \omega f, \theta h_i \right\rangle  = 0.$$
		Conversely, let $\left\langle  f, \Tilde{k}_{\theta,i}\right\rangle  = 0$, then $\left\langle \omega f, \theta h_i \right\rangle =  0$, because $\omega f \in H^p_+$. Consequently, it remains to prove that $\left\langle \omega f, \theta g \right\rangle = 0$ for any $g\in H^q_+$, where $q=\frac{p}{p-1}$. Clearly, $H^q_+ =K^q_\omega + \omega H^q_+$, so by Lemma \ref{sp} we have $H^q_+ = \omega H^q_+ + \mathrm{span} \langle h_i \rangle$. Therefore $g =\omega h +  h_i$ for some $h \in H^q_+$, so that entails
		$$\left\langle  \omega f, \theta g \right\rangle = \left\langle  f, \theta h \right\rangle + \left\langle \omega f, \theta h_i \right\rangle  = 0.$$
		This concludes the statement.
	\end{proof}
	
	For simplicity, throughout this paper we fix notation for the conjugate kernel $\tilde{k}_{\theta,i}$, namely denote $\mathsf{k} = \tilde{k}_{\theta,i}$. Also, let $(\mathrm{span}\langle\mathsf{k} \rangle )^\perp$ be the orthogonal complement of $\mathrm{span}\langle\mathsf{k} \rangle$ in the space $K_\theta^q$, so $(\mathrm{span}\langle\mathsf{k} \rangle )^\perp$ is a subspace of $K_\theta^p$. 
	
	Define the projector $\mathbb{K} \,\colon\, K_\theta^p \to (\mathrm{span}\langle\mathsf{k} \rangle )^\perp$ by
	$$\mathbb{K}\,\colon\, f\mapsto f -\left\| \mathsf{k} \right\|^{-2}_{L^2(\R)}\, \left\langle f,\mathsf{k} \right \rangle \mathsf{k}.$$
	The usual tensor notation will be used for operators of rank one: $f\otimes g$ denotes the operator defined by $(f\otimes g)[h]=\langle h,g\rangle f$. So, one can rewrite $\mathbb{K}= I - a\, \mathsf{k} \otimes \mathsf{k}$ with $I$ being the identity map on $K_\theta^p$ and $a = \left\| \mathsf{k} \right\|^{-2}_{L^2(\R)}$. Remark that using this new notation and Proposition \ref{l3}, the condition $f,\omega f\in K_\theta^p$ from Corollary \ref{con3} turns into $f \in \mathrm{Ran}\, \mathbb{K}$. Due to duality, we also get that $g,\omega g\in K_\theta^q$ turns into $f \in \mathrm{Ran}\, \mathbb{K}^*$.

	\medskip

	\subsection{Toeplitz commutator theorem}
	Let $k_{\theta,z}^\mathbb{D},\,z\in\D,$ be a reproducing kernel in the space $K_\theta^2(\mathbb{D})$, see e.g., \cite{saras}. For $z\in\mathbb{D}$ define the inverse conformal map:
	$$\omega^{-1} \,\colon\,z\mapsto i \cdot\frac{1+z}{1-z}\in\mathbb{C}_+.$$ 
	Note that $\omega^{-1}$ is called M{\"o}bius transformation and $\omega^{-1}(z)|_{z\in\T\setminus\{1\}}\in\R$.
	Now, for the operator $\mathbb{U}$ (see \ref{prem1}), its inverse map from $L^p(\R)$ to $L^p(\T)$ can be defined by
	$$\mathbb{U}^{-1}[f](z) = \left(\frac{-4\pi}{(1-z)^2}\right)^{\frac{1}{p}} f(\omega^{-1}(z)), \qquad z \in \T\setminus\{1\}, \, f\in L^p(\R).$$
	Further in the paper, we fix the following notation:
	$$\eta = a\cdot(-\bar{h}_i)^{\frac{2}{p}}\in L^p(\R),\quad p\in [1,+\infty).$$ 
	
	We prove auxiliary lemmas and, finally, conclude Corollary \ref{con3}.
	 
	\begin{Lemma}\label{tcker} \textit{
		The following statements hold true.
		\begin{itemize}
			\item[(a)] $k_{\theta, i}^{\C_+} = (-4\pi)^{-\frac{1}{p}}\, \mathbb{U}[k_{\theta,0}^\D]$.
			\item[(b)] Denote $\rho(z) = \pi^{\frac{1}{p}}\sigma(z),\,z\in\C_+$, and let $A_\phi$ be a bounded Toeplitz operator on $K_\theta^p(\D)$ with symbol $\phi \in K_\theta^p(\mathbb{D})$. Then, $T_{\rho\, \mathbb{U}[\phi]} =\mathbb{U}\, A_\phi \,\mathbb{U}^{-1}$.
			\item[(c)] For any $\varphi \in \sigma K^p_\theta$, we have $T_{\bar{\varphi}}[\mathsf{k}] =\frac{1}{a} \eta\cdot\tilde{\varphi}$.
		\end{itemize}}
	\end{Lemma}
	\begin{proof} We prove (a) and (b) independently and then prove (c) by the first two~statements.
	\begin{itemize}
		\item[(a)]
		Let $f \in L^p(\mathbb{D})$, then we get
		$$\mathbb{U}^{-1}\mathbb{U}[f](z)=f(z) = \left\langle f, k_{\theta, z}^\D \right\rangle =  \left\langle \mathbb{U}[f], \mathbb{U}[k_{\theta, z}^\D] \right\rangle.$$
		Also, $\mathbb{U}[f](i) =(-4\pi)^{-\frac{1}{p}}f(0)$. Therefore, 
		$$\left\langle \mathbb{U}[f], k_{\theta, i}^{\C_+} \right\rangle = \mathbb{U}[f](i) = (-4\pi)^{-\frac{1}{p}} \left\langle \mathbb{U}[f], \mathbb{U}[k_{\theta, 0}^\D] \right\rangle.$$
		
		\item[(b)] We note that the following is true almost everywhere in $\R$,
		$$\mathbb{U}[\phi\, \mathbb{U}^{-1} [f]] = \left(\frac{-4\pi}{(1-\omega)^2}\right)^{\frac{1}{p}}\mathbb{U}[\phi]\, f = \rho\, \mathbb{U}[\phi] \, f, $$
		since for any $x\in\R$,
		$$(-4\pi)^{\frac{1}{p}}\left(\frac{1}{1-\frac{x-i}{x+i}}\right)^{\frac{2}{p}}=(-4\pi)^{\frac{1}{p}}\left(\frac{x+i}{2i}\right)^{\frac{2}{p}}=\pi^{\frac{1}{p}}\sigma(x).$$
		Thus, for any $f\in K^p_\theta, g\in K^q_\theta$ we get
		$$\left\langle \mathbb{U}\, A_\phi\, \mathbb{U}^{-1} [f], g \right\rangle = \left\langle \mathbb{U}[\phi\, \mathbb{U}^{-1} [f]], g \right\rangle = \left\langle \rho\, \mathbb{U}[\phi] \, f, g \right\rangle = \left\langle T_{\rho\, \mathbb{U}[\phi]}[f], g \right\rangle.$$
		
		\item[(c)] By the properties of $\mathrm{P}_\theta$ on $K_\theta^p(\D)$ we obtain 
		$$A_\phi[k_{\theta, 0}^\D] =\mathrm{P}_{\theta}[ \phi (1-\overline{{\theta}(0)}\theta)] = \phi - \overline{\theta (0)}\cdot\mathrm{P}_{\theta}[\phi \theta] = \phi.$$
		Summarizing all of the above, for almost every $x\in\R$,
		\begin{align*}
			T_{\varphi}[k_{\theta, i}^{\C_+}](x) &= \mathbb{U}A_{\mathbb{U}^{-1}[\frac{\varphi}{\rho}]}\mathbb{U}^{-1}[k_{\theta,0}^{\C_+}](x)\\ &=(-4\pi)^{-\frac{1}{p}}\, \mathbb{U}A_{\mathbb{U}^{-1}[\frac{\varphi}{\rho}]} [k_{\theta,0}^{\D}](x)\\
			&= (-4\pi)^{-\frac{1}{p}}\, \mathbb{U}\mathbb{U}^{-1}\left[\frac{\varphi}{\rho}\right](x) \\
			&=(-4)^{-\frac{1}{p}}(\pi(x+i))^{-\frac{2}{p}}\,\varphi(x) .
		\end{align*}
		The space $K_\theta^p$ is closed under the conjugation $C\,\colon\,f\mapsto \theta \bar{f} = \tilde{f}$. Recall that by the definition of the conjugate kernel we have $C[\mathsf{k}]=k_{\theta,i}^{\C_+}$. Thus by conjugate-symmetric property (see \cite{saras}), i.e. $CT_\varphi C = T_{\bar{\varphi}}$, we obtain
		$$T_{\bar{\varphi}}[\mathsf{k}](x) = C T_\varphi  [{k}_{\theta, i}^{\C_+}](x) = (-2\pi i (x-i))^{-\frac{2}{p}} \, \tilde{\varphi}(x)$$
		for almost every $x\in\mathbb{R}$. This implies the required formula.
	\end{itemize}
	\end{proof}
	
	Denote two particular Toeplitz operators $T_{\omega}$ and $T_{\bar{\omega}}$ as $\Lambda$ and $\bar{\Lambda}$ respectively, i.e., $\Lambda = T_{\omega}$ and $\bar{\Lambda} = T_{\bar{\omega}}$. So $\Lambda$ is a continuous case version of the truncated Toeplitz operator $T_z$ in the discrete case, which is the compression of the forward shift operator $S[f](z)=zf(z)$ from $H^p(\D)$ to $K_\theta^p(\D)$. Furthermore, $\bar{\Lambda}$ is a continuous case version of the truncated Toeplitz operator $T_{\bar{z}}$ in the discrete case, which is the compression of the backward shift operator $S^*[f](z)=\frac{1}{z}(f(z)-f(0))$ to $K_\theta^p(\D)$, see e.g., \cite{cima}, \cite{garnet}, \cite{nikk}.
	
	\begin{Lemma}\label{ktimek}
		\textit{We have $I-\bar{\Lambda} \Lambda =a\,\mathsf{k} \otimes \mathsf{k}$ on $K_\theta^p$.}
	\end{Lemma}
	\begin{proof}
		Indeed, for any $f \in K^p_\theta$ and $g \in \mathrm{Ran}\,\mathbb{K}^{*}$ we have $$\left\langle \bar{\Lambda}\Lambda [f], g \right\rangle = \left\langle \mathbb{P}_\theta[\omega f], \mathbb{P}_\theta^* [\omega g] \right\rangle = \left\langle f, g \right\rangle.$$
		This implies that $\mathbb{K}(I-\bar{\Lambda} \Lambda)=0$ on $K_\theta^p$. Hence $\mathrm{Ran}\,(I -\bar{\Lambda} \Lambda)  \subset \mathrm{span}\langle \mathsf{k}\rangle$, and then we obtain $I-\bar{\Lambda} \Lambda =a \,\mathsf{k} \otimes \mathsf{k}$. 
	\end{proof}
	
	\begin{Lemma}\label{comm} \textit{
		For any $\phi \in \sigma K_\theta^p$ and $\psi \in \sigma K_\theta^q$, we have
		$$T_{\bar{\phi}+\psi} - \bar{\Lambda} T_{\bar{\phi}+\psi} \Lambda= \eta \Tilde{\phi} \otimes \mathsf{k} + \mathsf{k} \otimes \eta \Tilde{\psi} .$$}
	\end{Lemma}
	\begin{proof}
		We know that $H^p_+ =K^p_\theta + \theta H^p_+$, then for any function $h \in H^p_+$ there exists $g \in H^p_+$ such that $\mathbb{P}_\theta[h] = h + \theta g$. Hence, for any $f\in K_\theta^p$ and some $g_1,g_2 \in H^p_+$ we have
		$$T_\psi T_\omega [f] = T_\psi [\omega f + \theta g_1] = T_{\psi \omega}[f] =T_\omega [\psi f + \theta g_2] = T_\omega T_\psi [f].$$
		Therefore, $T_{\psi} \Lambda = \Lambda T_{\psi}$ and by duality $\bar{\Lambda} T_{\bar{\phi}} = T_{\bar{\phi}} \bar{\Lambda}$. Then, Lemma \ref{ktimek} entails
		\begin{align*}
			T_{\bar{\phi}+\psi} - \bar{\Lambda} T_{\bar{\phi}+\psi} \Lambda &=   T_{\bar{\phi}} (I-\bar{\Lambda}\Lambda) + (I-\bar{\Lambda}\Lambda) T_{\psi}\\
			&=a\,T_{\bar{\phi}}\,( \mathsf{k} \otimes \mathsf{k}) + a \,(\mathsf{k} \otimes \mathsf{k}) T_\psi.
		\end{align*}
		Note that by duality $(\mathsf{k} \otimes \mathsf{k}) T_\psi = \mathsf{k} \otimes T_{\bar{\psi}}[\mathsf{k}]$. 
		Thus, Lemma \ref{tcker} yields the required identity.
	\end{proof}
	
	Recall that the space of all bounded Toeplitz operators on $K_\theta^p$ is denoted by $\mathcal{T}^p(\theta)$. 
 
	\begin{Proposition}\label{ex_sym} \textit{
		Given a bounded operator $T \,\colon\, K_\theta^p\to K_\theta^p$. Suppose there are functions $\phi \in \sigma K_\theta^p,\,\psi \in \sigma K_\theta^q$ such that 
		$$T - \bar{\Lambda} T \Lambda = \eta \Tilde{\phi} \otimes \mathsf{k} + \mathsf{k} \otimes \eta \Tilde{\psi},$$
		then $T$ belongs to $\mathcal{T}^p(\theta)$, in which case $T = T_{\bar{\phi}+\psi}$.}
	\end{Proposition} 
	\begin{proof}
		Let $\mathcal{C}_{\bar{\phi}+\psi} = T_{\bar{\phi}+\psi} - \bar{\Lambda} T_{\bar{\phi}+\psi} \Lambda$ be a Toeplitz commutator. Then for every integer $N \geqslant 0$ and any $f\in K_\theta^p$, $g\in K_\theta^q$, we have
		$$\left\langle T_{\bar{\phi}+\psi}[f], g\right\rangle =\sum\limits_{n=0}^{N} \left\langle\bar{\Lambda}^n\mathcal{C}_{\bar{\phi}+\psi} \Lambda^n[f], g \right\rangle + \left\langle T_{\bar{\phi}+\psi} \Lambda^{N+1} [f], \Lambda^{N+1}[g]\right\rangle.$$ 
		By the commutation property of $T_\psi$ and $\Lambda$ from the proof of Lemma \ref{comm}, the last term on the right side can be presented as 
		$$\left\langle \Lambda^{N+1}T_{\psi}[f], \Lambda^{N+1}[g]\right\rangle + \left\langle \Lambda^{N+1}[f], \Lambda^{N+1}T_{\phi}[g]\right\rangle.$$
		It is straightforward to show that $\Lambda^N \to 0$ in the weak operator topology and $\bar{\Lambda}^N \to 0$ in the strong operator topology as $N\to +\infty$ (since for the shift operators we have $S^N\to 0$ in the weak operator topology and $(S^*)^N\to 0$ in the strong operator topology). This entails that $T_{\bar{\phi}+\psi} = \sum_{n\geqslant 0} \bar{\Lambda}^n\mathcal{C}_{\bar{\phi}+\psi} \Lambda^n$ and the series converges in the strong operator topology. 
		
		It remains to prove that operator $T$ can be represented in the same form. Let $\mathcal{C} = T - \bar{\Lambda}T \Lambda$ be a commutator of $T$ and note that $\mathcal{C} = \mathcal{C}_{\bar{\phi}+\psi}$ by the condition from the proposition and the previous Lemma \ref{comm}. For any integer $N \geqslant 0$ we have 
		$$T = \sum_{n=0}^N \bar{\Lambda}^n\mathcal{C} \Lambda^n + \bar{\Lambda}^{N+1} T \Lambda^{N+1} =  \sum_{n=0}^N \bar{\Lambda}^n\mathcal{C}_{\bar{\phi}+\psi} \Lambda^n + \bar{\Lambda}^{N+1} T \Lambda^{N+1}.$$
		The latter summand on the right side 
		tends to $0$ as $N\to\infty$, again due to the convergence $\bar{\Lambda}^N \to 0$ in the strong operator topology. Thus, $T=T_{\bar{\phi}+\psi}$ and the proposition is proved. 
	\end{proof}
	
	Now, we are ready to prove that the condition $\left\langle T[f], g \right\rangle = \left\langle T [\omega f],  \omega g\right\rangle$ for any $f,\omega f\in K_\theta^p,\,g,\omega g\in K_\theta^q$, in Corollary \ref{con3} is sufficient for $T$ to be a Toeplitz operator on $K_\theta^p$. The necessity is immediate.
	
	\begin{Theorem}\label{propertyT} \textit{
		For any bounded operator $T \,\colon\, K_\theta^p\to K_\theta^p$ such that 
		$$\left\langle T[f], g \right\rangle = \left\langle T [\omega f],  \omega g\right\rangle, \qquad  f \in \mathrm{Ran}\,\mathbb{K}, g \in \mathrm{Ran}\,\mathbb{K}^{*},$$
		there exists a symbol $\varphi \in (\overline{\sigma K^p_\theta}+\sigma K^q_\theta)$ such that $T = T_\varphi\in\mathcal{T}^p(\theta)$.}
	\end{Theorem}
	\begin{proof}
		Let $\mathcal{C} = T - \bar{\Lambda} T \Lambda$. Notice that $\left\langle \mathcal{C}[f], g\right\rangle = \left\langle T[f], g\right\rangle- \left\langle T [\omega f],  \omega g\right\rangle =0$, then $\mathbb{K}\,\mathcal{C}\,\mathbb{K} = 0$ on $K_\theta^p$, which implies $(I-a \mathsf{k}\otimes \mathsf{k})(\mathcal{C}-a \mathsf{k}\otimes \mathsf{k})=0$. So, we obtain 
		$$\mathcal{C} = a\,\mathcal{C}[\mathsf{k}] \otimes \mathsf{k}  +  a\,\mathsf{k}\otimes \mathcal{C}^*[\mathsf{k}] - a^2\,\left\langle \mathcal{C}[\mathsf{k}], \mathsf{k} \right \rangle\mathsf{k} \otimes \mathsf{k}.$$
		Denote $Q = \mathcal{C}^* - \bar{a}\,\overline{\left\langle \mathcal{C}[\mathsf{k}], \mathsf{k} \right \rangle} I$. This operator acts from $K_\theta^q$ to $K_\theta^q$. Thus, one can rewrite the identity above:
		$$\mathcal{C} =a\, \mathcal{C}[\mathsf{k}]\otimes \mathsf{k} +  a\,\mathsf{k}\otimes Q[\mathsf{k}].$$
		Take $\phi \in \sigma K_\theta^p$ and $\psi \in \sigma K_\theta^q$ such that $\eta\tilde{\phi} = a\,\mathcal{C}[\mathsf{k}]$ and $\eta\tilde{\psi} = a\,Q[\mathsf{k}]$, hence Proposition~\ref{ex_sym} proves the theorem.
	\end{proof}
	
	A result similar to Theorem \ref{propertyT} is usually called the Toeplitz commutator theorem.	
	
	\medskip
	
	\section{Duality methods and weak factorization theorem}\label{7}
	Let $1<p<+\infty$ and $q$ be the H{\"o}lder conjugate of $p$. Here we again consider only the class of bounded Toeplitz operators from $\mathcal{T}^p(a)$, symbols are taken from $\mathcal{P}(\R)$. We recall that $\pw=\bar{\theta}_aK_{\theta^2_a}^p$ and $\P_a=\bar{\theta}_a\P_{\theta_a^2}$, where $\theta_a(z)=e^{2\pi i a z}$. 
	
	A special version of Theorem \ref{propertyT} is as follows.
	
	\begin{Theorem}\label{keypw} \textit{
			A bounded operator $T\,:\,\pw\to\pw$ belongs to $\mathcal{T}^p(a)$ if and only if the condition $f,\omega f \in \pw,\, g,\omega g \in \mathrm{PW}_a^q$ yields
			$$\left\langle T[f], g \right\rangle = \left\langle T [\omega f],  \omega g\right\rangle.$$}
	\end{Theorem}
	
	Define the following special subspace of $\mathrm{PW}_{2a}^1$:
	$$\mathcal{X}^{p,q} = \left\{\sum_{k\geqslant 0} {f_{k} \bar{g}_{k}}\; \Biggm\vert\; f_{k} \in \pw, \, g_k \in \mathrm{PW}_a^q,\, \sum_{k\geqslant 0} \left \| f_k \right \|_{L^p(\mathbb{R})} \cdot\left \| g_k \right \|_{L^q(\mathbb{R})} <+\infty \right\}.$$
	The norm in $\mathcal{X}^{p,q}$ is defined as the infimum of $\sum_{k\geqslant 0} \left \| f_k \right \|_{L^p(\mathbb{R})} \cdot\left \| g_k \right \|_{L^q(\mathbb{R})}$ over all representations of the element in the form $\sum_{k\geqslant 0} f_k \bar{g}_k$. This norm makes $\mathcal{X}^{p,q}$ a Banach space. We show that $\mathcal{X}^{p,q}$ is a predual to the space of all bounded Toeplitz operators:
	$$(\mathcal{X}^{p,q})^* \cong \mathcal{T}^p(a),$$
	where the notion $\cong$ means that the spaces above are isometrically isomorphic. 
	
	\begin{Proposition}\label{dualX} \textit{
		The dual space $(\mathcal{X}^{p,q})^*$ can be naturally identified with $\mathcal{T}^p(a)$ and all continuous linear functionals over $\mathcal{X}^{p,q}$ are of the form:
		$$\Phi(h) = \sum_{k\geqslant 0} \left \langle  T[f_k] , g_k \right \rangle, \qquad h = \sum_{k\geqslant 0} f_k \bar{g}_k \in \mathcal{X}^{p,q}, $$
		with $T \in \mathcal{T}^p(a)$, and the correspondence between the functionals from $(\mathcal{X}^{p,q})^*$ and the space $\mathcal{T}^p(a)$ is a one-to-one isometry.}
	\end{Proposition}
	\begin{proof}
		First, we verify that the functional is well defined for an operator $T_\varphi \in \mathcal{T}^p(a)$, that is, the value of a functional is independent of the particular representation chosen for $h\in\mathcal{X}^{p,q}$. Suppose $h =\sum_k f_k \bar{g}_k= 0$, then
		$$\Phi(h) = \int\limits_{\R}  \varphi\left(\sum_{k\geqslant 0} f_k \bar{g}_k\right) \, d\lambda = \int\limits_{\R} \varphi  h \,d\lambda =0,\qquad \varphi\in\mathcal{P}(\R),$$
		where $\lambda$ is the Lebesgue measure on $\R$.
		
		Now, prove the equality $\left \|T \right\|_{\pw\to \pw} = \left \|\Phi \right\|$. Indeed, by the definition of the functional it is obvious that $\left \|\Phi \right\|\leqslant \left \|T \right\|$. Also, for any unit norm vectors $f \in \pw$ and $g \in \mathrm{PW}_a^q$ we have $\left \| f\bar{g} \right \|_{\mathcal{X}^{p,q}} \leqslant 1$ and 
		$$\left \| T \right \| = \sup\limits_{\left \| f \right \|_{L^p(\R)},\left \| g \right \|_{L^q(\R)} \leqslant 1} |\left \langle T[f], g \right \rangle| = \sup\limits_{\left \| f \right \|_{L^p(\R)},\left \| g \right \|_{L^q(\R)}\leqslant 1} | \Phi( f \bar{g} )| \leqslant \left \| \Phi \right \|.$$
		This proves the inverse inequality.
		
		It remains to show that any linear continuous functional $\Phi \in (\mathcal{X}^{p,q})^*$ may be represented in the form $\Phi = \sum_k \left \langle  T \cdot , \cdot \right \rangle$ for some unique Toeplitz operator $T \in \mathcal{T}^p(a)$.
		Pick a continuous functional $\Phi\in (\mathcal{X}^{p,q})^*$ and define the operator $T\, \colon\, \pw\to\pw $ by its sesquilinear form $\left \langle  T [f], g \right \rangle =  \Phi(f\bar{g} )$ for any $f \in \pw$ and $g \in \mathrm{PW}_a^q$. If $\omega f\in \pw$ and $\omega g \in \mathrm{PW}_a^q$, then 
		$$\left \langle  T [f], g \right \rangle =  \Phi(f\bar{g} ) = \Phi(\omega f\bar{\omega} \bar{g} ) = \left \langle  T [\omega f], \omega g \right \rangle.$$
		Thus, by Theorem \ref{keypw}, we obtain that $T \in \mathcal{T}^p(a)$. The uniqueness of $T$ is a consequence of the relation $\left \|T \right\| = \left \|\Phi \right\|$.
	\end{proof}
	
	\medskip
	
	\subsection{Schwartz symbol Toeplitz operators describe almost all bounded Toeplitz operators}\label{weakstar}
	Let $1<p<+\infty$ and $q$ be the H{\"o}lder conjugate of $p$. By Proposition \ref{dualX}, one can view any $T \in \mathcal{T}^p(a)$ as a bounded linear functional on $\mathcal{X}^{p,q}$, in which case, we write the dual paring as 
	$$\langle T,h\rangle = \sum_{k\geqslant 0}\langle T[f_k], g_k\rangle$$ 
	for every element $h=\sum_k f_k\bar{g}_k\in\mathcal{X}^{p,q}$. We equip $\mathcal{T}^p(a)$ with the weak${}^*$ topology, so continuous linear functionals are only those from $\mathcal{X}^{p,q}$ being treated as elements of $(\mathcal{T}^p(a))^*$ due to the canonical embedding  $\mathcal{X}^{p,q}\hookrightarrow ((\mathcal{X}^{p,q})^*)^*\cong  (\mathcal{T}^p(a))^*$. A Toeplitz operator with a symbol from the Schwartz space $\sr$ is called a Schwartz Toeplitz operator.
	
	\begin{Proposition}\label{weak*dense} \textit{
		Schwartz Toeplitz operators on $\pw$ are weak${}^*$ dense in~$\mathcal{T}^p(a)$.}
	\end{Proposition}
	\begin{proof}
		Clearly, the space of all Schwartz Toeplitz operators on $\pw$ separates the points of $\mathcal{X}^{p,q}$, since for an element $h\in\mathcal{X}^{p,q}$,
		$$0=\langle T_\varphi ,h\rangle = \langle \varphi, \bar{h}\rangle,\qquad \varphi \in\S(\R),$$
		and the fundamental lemma of the calculus of variations yield that $h=0$ in $\mathcal{X}^{p,q}$.
		
		To show the claim, let us suppose the opposite, $\{T_\varphi\in\mathcal{T}^p(a)\mid\varphi\in\S(\R)\}$ is not weak${}^*$ dense in $\mathcal{T}^p(a)$. Let $\mathcal{T}_\mathcal{S}$ be the weak${}^*$ closure of all Schwartz Toeplitz operators. Choose some $A\in\mathcal{T}^p(a)\setminus \mathcal{T}_\mathcal{S}$. Then, by the Hahn–Banach separation theorem, there exists an element $h\in\mathcal{X}^{p,q}$ being treated as a continuous linear functional $J_h\,:\, \mathcal{T}^p(a)\to\C$,
		$$J_h\,:\,T\mapsto \langle T, h\rangle,\qquad T\in\mathcal{T}^p(a) ,$$
		such that $J_h|_{\mathcal{T}_\mathcal{S}}=0$ and $J_h(A) \ne 0$. The result $\langle\cdot,h\rangle|_{\mathcal{T}_\mathcal{S}}=J_h|_{\mathcal{T}_\mathcal{S}} =0$ implies $h=0$, which immediately gives a contradiction, $0\ne J_h(A) =\langle A, h\rangle =0$. 
		\end{proof}
		
		Notice that Proposition \ref{dualX} implies
		\begin{equation}\label{eqhxpq}
			\left \| h \right \|_{\mathcal{X}^{p,q}} = \sup \left\{| \langle T,h\rangle| \;:\; T \in \mathcal{T}^p(a),  \| T \| \leqslant 1 \right\}.
		\end{equation}
		We write $a\asymp b$, if there are some constants $C_1,C_2>0$ that $C_1 b \leqslant a \leqslant C_2 b$.
			
		\begin{Corollary}\label{corasymp} \textit{The following norm equivalence holds:
			$$\left \| h \right \|_{\mathcal{X}^{p,q}} \asymp \sup \left\{| \langle T_\varphi,h\rangle| \;:\; T_\varphi \in \mathcal{T}^p(a),\varphi\in \S(\R),  \| T_\varphi \| \leqslant 1 \right\}.$$}
		\end{Corollary}
		\begin{proof}		
		Due to identity (\ref{eqhxpq}), it remains to show that for a universal $C>0$,
		\begin{equation}\label{rggds}
			\left \| h \right \|_{\mathcal{X}^{p,q}}\leqslant C\,\sup \left\{| \langle T_\varphi,h\rangle| \;:\; \varphi\in\S(\R),  \| T_\varphi \| \leqslant 1 \right\}.
		\end{equation}
		
		Proposition \ref{weak*dense} yields that for any $T\in\mathcal{T}^p(a)$, $\|T\|\leqslant 1$, and any $\varepsilon>0$, there is a family $F$ of Schwartz Toeplitz operators $T_\phi$ such that $$\left|\langle T, h\rangle -\langle T_\phi, h\rangle \right|<\varepsilon,\qquad h\in \mathcal{X}^{p,q}.$$ 
		Since the values $|\langle T_\phi, h\rangle|$ for the members $T_\phi$ of the family $F$ are pointwise bounded,
		$\left|\langle T_\phi, h\rangle \right| < \left|\langle T, h\rangle \right| +\varepsilon$ for any $h\in \mathcal{X}^{p,q}$, then the Banach–Steinhaus theorem implies the uniform boundedness, that is, for a constant $C\geqslant 1$,
		$$\sup\left\{\| T_\phi\| \,:\, T_\phi\in F\right\}\leqslant C.$$
		
		Fix a member $T_\phi$ of the family $F$, and note that also we have
		$$\left|\langle T, h\rangle \right| < \left|\langle T_\phi, h\rangle \right| +\varepsilon,\qquad h\in \mathcal{X}^{p,q}.$$
		Denote $T_{\phi/C} = \frac{1}{C}T_\phi$, so the inequality above becomes
		$$\frac{1}{C}\left|\langle T, h\rangle \right| < \left|\langle T_{\phi/C}, h\rangle \right| +\frac{\varepsilon}{C},\qquad h\in \mathcal{X}^{p,q}.$$
		Remark that $\|T_{\phi/C}\|\leqslant 1$. Hence, this gives us
		$$\frac{1}{C}\left|\langle T, h\rangle \right| < \left|\langle T_{\phi/C}, h\rangle \right| +\frac{\varepsilon}{C} \leqslant \sup \left\{| \langle T_\varphi,h\rangle| \;:\; \varphi\in\S(\R),  \| T_\varphi \| \leqslant 1 \right\} + \frac{\varepsilon}{C}.$$
		Multiply by $C$, then setting $\varepsilon\to 0$ entails
		$$\left|\langle T, h\rangle \right|\leqslant C\, \sup \left\{| \langle T_\varphi,h\rangle| \;:\; \varphi\in\S(\R),  \| T_\varphi \| \leqslant 1 \right\}. $$
		Finally, the arbitrariness of $T\in\mathcal{T}^p(a)$, $\|T\|\leqslant 1$, yields estimate (\ref{rggds}), which, in turn, concludes the Corollary.
	\end{proof}
	
	\medskip
	
	\subsection{Proof of the weak factorization theorem}	
	Finally, we prove the second result of the paper, namely the weak factorization Theorem \ref{result2}. Recall that we want to show that for any function $h\in \mathrm{PW}_{2a}^1$ there exist $f_k \in \mathrm{PW}_a^p$, $g_k \in \mathrm{PW}_a^q$ with $$\sum_{k\geqslant 0} \left \|f_k\right \|_{L^p(\R)}\cdot\left \|g_k\right \|_{L^q(\R)}<+\infty \;\text{ such that }\;h = \sum_{k\geqslant 0} f_k g_k,$$
	where $1 < p < \infty$ and $\tfrac{1}{p}+\tfrac{1}{q}=1$. 
	
	This statement is equivalent to the claim that $\mathrm{PW}_{2a}^1 = \mathcal{X}^{p,q}$. In the general case, we have $\mathrm{PW}_{2a}^1=\cl_{L^1(\R)} \mathcal{X}^{p,q}$.	So, $\mathcal{X}^{p,q}$ must be a complete space in $L^1(\R)$-norm. However, we proceed with a much simpler proof. 
	
	\begin{proof}
		
		By the Closed Graph theorem, $\mathrm{PW}_{2a}^1= \mathcal{X}^{p,q}$ if and only if the norms in $\mathcal{X}^{p,q}$ and $L^1(\R)$ are equivalent. The embedding $\mathcal{X}^{p,q} \hookrightarrow \mathrm{PW}_{2a}^1$ is clearly continuous, $$\left\| h\right\|_{L^1(\mathbb{R})} \leqslant \left\| h\right\|_{\mathcal{X}^{p,q}},\qquad h \in \mathcal{X}^{p,q}.$$ 
		Therefore, we are to prove that the inverse embedding is continuous as well.
		
		Corollary \ref{corasymp} implies that the set of all Schwartz Toeplitz operators is sufficient to describe the norm on $\mathcal{X}^{p,q}$. Theorem \ref{result} guarantees that for any Schwartz Toeplitz operator $T_\varphi$ on $\pw$, there is a bounded symbol $\psi\in L^\infty(\R)$ and a constant $c_p>0$ such that $T_\varphi=T_\psi$ and
		$$\left \| \psi \right\|_{L^\infty(\mathbb{R})} \leqslant c_p \left \|T_\psi \right \|.$$ Consequently, for every $h \in \mathcal{X}^{p,q}$ we obtain the required upper bound:
		\begin{align*}
			\left \| h \right \|_{\mathcal{X}^{p,q}} 
			&\leqslant C\,\sup \left\{\left| \left  \langle  T_\varphi, h \right\rangle\right| \;:\; \varphi\in\S(\R),  \| T_\varphi \| \leqslant 1 \right\}\\
			&\leqslant C\,\sup \left\{\left|  \left  \langle  T_\psi, h \right\rangle\right| \;:\; \psi\in L^\infty(\R),  \| T_\psi \| \leqslant 1 \right\}\\
			&= C\,\sup \left\{\left|  \left  \langle  \psi , \bar{h}  \right\rangle\right|\;:\;  \left \| \psi \right\|_{L^\infty(\R)} \leqslant c_p \right\}\\
			&= C_p\, \sup\left\{ \left|  \left  \langle  \psi, \bar{h}  \right\rangle\right|\,:\,\left \| \psi \right\|_{L^\infty(\mathbb{R})} \leqslant 1\right\}\\
			& \leqslant C_p\,\sup\limits_{\left \| \psi \right\|_{L^\infty(\mathbb{R})} \leqslant 1} \;\int\limits_{\R} |\psi h| \,d\lambda \;\; \leqslant\;\; C_p\,\left \| h \right \|_{L^1(\mathbb{R})},
		\end{align*}
	where $C_p=C\cdot c_p$ and $\lambda$ is the Lebesgue measure on $\R$. This proves the theorem.
	\end{proof}
	
	\appendix
	\section*{Appendix}\label{appendix}
	We give a sketch of the proof of Theorem \ref{nehari}. Recall the statement.
	
	Let  $1<p<+\infty$ and let $\varphi \in L^\infty(\R)$. Then there exists a function $\psi\in L^\infty(\R)$ such that $H_\psi=H_\varphi$ and, moreover, $\left \| \psi \right \|_{L^\infty(\R)} \leqslant \left \| H_\varphi \right \|_{H^p_+\to H^p_-}$.
	\begin{proof}
		Consider a function $\varphi \in L^\infty(\R)$. We have 
		$$\left \| H_\varphi \right \| = \sup \{\langle \varphi f , \P_-[g] \rangle \mid f\in H^p_+, g\in L^q(\R),\;\left \| f \right \|_{L^p(\R)} \leqslant 1,\left \| g\right \|_{L^q(\R)} \leqslant 1 \},$$
		where $\tfrac{1}{p}+\tfrac{1}{q}=1$. Choosing $g\in H^q_-$ we see that
		$$\left \| H_\varphi \right \| \geqslant \sup \{\langle \varphi, \overline{fh} \rangle \mid f\in H^p_+, h\in H^q_+ ,\; \left \| f \right \|_{L^p(\R)} \leqslant 1,\left \| h\right \|_{L^q(\R)} \leqslant 1 \}.$$
		Since every function $F$ in the unit ball of $H^1_+$ can be represented in the form $F = f h$ for some $f\in H^p_+, h\in H^q_+$, we have 
		$$\left \| H_\varphi \right \| \geqslant \sup \{\langle \varphi, \overline{F} \rangle \mid F\in H^1_+ ,\; \left \| F \right \|_{L^1(\R)} \leqslant 1 \}.$$
		Extending the linear functional $\Phi_\varphi \, \colon \, F \to \langle \varphi, \overline{F}\rangle$ from $H^1_+$ to $L^1(\R)$ by the Hahn-Banach theorem, we see that there exists a function $\psi \in L^\infty(\R)$ such that $\left \| \psi \right \|_{L^\infty(\R)} \leqslant \left \| H_\varphi \right \|$ and $\langle \varphi, \overline{F}\rangle = \langle \psi, \overline{F} \rangle$ for every $F \in H^1_+$. In particular, we have $\langle \varphi f, g \rangle = \langle \psi f, g \rangle$ for all $f\in H^p_+$, $g\in H^q_-$. In other words, \text{$H_\varphi = H_\psi$.}
	\end{proof}
	
	\medskip
	
	\begin{center}
		{\sc Acknowledgement}
	\end{center}
	I am grateful to R. Bessonov for constant attention to this work, various helpful suggestions and discussions.
	
	This work was partially supported by the Finnish Ministry of Education and Culture’s Pilot for Doctoral Programmes (Pilot project Mathematics of Sensing, Imaging and Modelling), the FAME flagship of the Research Council of Finland (Flagship of Advanced Mathematics for Sensing Imaging and Modelling, grant 359182), the Centre of Excellence of Inverse Modelling and Imaging of Research Council of Finland (grant 336786), and the European Research Council (the Advanced Grant project
	101097198). Views and opinions expressed are those of the author only and do not necessarily reflect	those of the European Union or the other funding organizations.
		
	\bibliographystyle{plain} 
	\bibliography{references}	
	
	\vspace{1.5cm}
	\begin{flushleft}
		\begin{footnotesize}
			\sc University of Helsinki
			
			Department of Mathematics and Statistics
			
			Pietari Kalmin katu 5, 00560, Helsinki, Finland
			
			\textit{E-mail address:} \tt{petr.kulikov@helsinki.fi}
		\end{footnotesize}
	\end{flushleft}
\end{document}